\documentclass[ejs]{imsart}
\RequirePackage[OT1]{fontenc}
\RequirePackage{amsthm,amsmath}
\RequirePackage[numbers]{natbib}
\RequirePackage[colorlinks,citecolor=blue,urlcolor=blue]{hyperref}

\pubyear{0}
\volume{0}
\issue{0}
\firstpage{1}
\lastpage{29}
\usepackage{graphicx,psfrag,epsf}
\usepackage{pdfsync}
\usepackage{amssymb}
\usepackage[latin1]{inputenc}
\usepackage{amsthm}

{\catcode `\@=11 \global\let\AddToReset=\@addtoreset}
\AddToReset{equation}{section}

\usepackage{makeidx}
\makeindex
\def\R{\mathbb{R}}

\def\N{\mathbb{N}}
\def\P{\mathbb{P}}
\def\E{\mathbb{E}}
\def\L{\mathbb{L}}

\def\R{\mathbb{R}}

\def\Z{\mathbb{Z}}

\def\1{\mbox{I\hspace{-.6em}1}} 
\def\v{\mbox{Var}\,}

\def\cov{\mbox{Cov}\,}

\def\Lip{\mbox{Lip}\,}
\def\1{\mbox{\hspace{.2em}I\hspace{-.6em}1}} 

\def\limiten{\renewcommand{\arraystretch}{0.5}
\begin{array}[t]{c}\stackrel{}{\longrightarrow} \\
{\scriptstyle
n\rightarrow+\infty}\end{array}\renewcommand{\arraystretch}{1}}

\def\limiteloin{\renewcommand{\arraystretch}{0.5}
\begin{array}[t]{c}\stackrel{{\cal L}}{\longrightarrow} \\
{\scriptstyle
n\rightarrow+\infty}\end{array}\renewcommand{\arraystretch}{1}}

\def\limiteproban{\renewcommand{\arraystretch}{0.5}
\begin{array}[t]{c}\stackrel{{\P}}{\longrightarrow} \\
{\scriptstyle
n\rightarrow+\infty}\end{array}\renewcommand{\arraystretch}{1}}

\newcommand{\findem}{\hfill\hbox{\hskip 4pt
\vrule width 5pt height 6pt depth 1.5pt}\vspace{.5cm}\par}

  \newtheorem{prop}{Proposition}[section]
  
\newtheorem{theo}{Theorem}[section]
 \newtheorem{lem}{Lemma}[section]
 
  \newtheorem{Rem}{Remark}[section]
\newtheorem{Def}{Definition}[section]

\begin{document}

\begin{frontmatter}
\title{Non-parametric estimation of time varying AR(1)--processes with local stationarity and periodicity}
\runtitle{Local periodic time varying AR(1)--processes}

\begin{aug}

\author{\fnms{Jean-Marc} \snm{Bardet}\thanksref{a,e1}\ead[label=e1,mark]{bardet@univ-paris1.fr}}
\and
\author{\fnms{Paul} \snm{Doukhan}\thanksref{b,e2}\ead[label=e2,mark]{doukhan@u-cergy.fr}}

\address[a]{SAMM EA4543, University Panth\'eon-Sorbonne, 90, rue de Tolbiac, 75634, Paris, France.
\printead{e1}}

\address[b]{AGM-UMR8088, University Cergy-Pontoise, France, and CIMFAV, Valparaiso, Chile.
\printead{e2}}

\runauthor{Bardet and Doukhan}

\affiliation{Some University and Another University}

\end{aug}

\begin{abstract}
Extending the ideas of \cite{dah}, this paper aims at providing a kernel based non-parametric estimation of a new class of time varying AR$(1)$ processes $(X_t)$, with local stationarity and periodic features (with a known period $T$), inducing the definition $X_t=a_t(t/nT)X_{t-1}+\xi_t$  for $t\in \N$ and with $a_{t+T}\equiv a_t$. Central limit theorems are established for kernel estimators $\widehat a_s(u)$  reaching classical minimax rates and only requiring low order moment conditions of the white noise $(\xi_t)_t$ up to the second order.
\end{abstract}



\end{frontmatter}

\

{\it This paper is dedicated to the memory of Jean Bretagnolle}

\section{Introduction}
Since the seminal paper \cite{dahl96}, the local-stationarity property provides new models and approaches for introducing non-stationarity in times series. The recently published handbook \cite{dah} gives a complete survey about new results obtained since $20$ years on this topics. \\
An interesting new kind of models is obtained from a natural extension of usual ARMA processes, so called 
tvARMA($p,q$)--processes  defined in \cite{dahl2009}, as:
\begin{equation}  \label{tvARMA} 
\sum_{j=0}^p \alpha_j\Big ( \frac t n \Big )\, X^{(n)}_{t-j}= \sum_{k=0}^q \beta_k \Big ( \frac t n \Big )\, \xi_{t-k},\qquad 1\le t\le n,
\end{equation} 
where $\alpha_j$ and $\beta_k$ are bounded functions. This is a special case of locally stationary linear process defined by $X^{(n)}_{t}=\sum_{j=0}^\infty \gamma_j\Big ( \frac t n \Big ) \, \xi_{t-j}$. Such models have been studied in many papers, especially concerning the parametric, semi-parametric or non-parametric estimations of functions $\alpha_j$, $\beta_k$ or $\gamma_j$, or other functions depending on these functions; see, for instance references \cite{dahl97}, \cite{dahl2009}, \cite{dah}, or \cite{grenier}, \cite{bf2003}, \cite{fg2004}, \cite{mpr2005} or \cite{am2006}. \\
~\\
For simplicity, we restrict in this first work  to time-varying AR$(1)$--processes $(X^{(n)}_t)$ including  a periodic component:
\begin{equation}\label{AR1}
X^{(n)}_t=a_{t}\left(\frac t{nT}\right)X^{(n)}_{t-1}+\xi_t,\quad \mbox{with $a_{t+T}\equiv a_t$},\ \mbox{for any} \left \{ \begin{array}{l} 1\le t\le nT \\ n\in \N^* \end{array} \right .,
\end{equation}
where $T\in \N^*$ is a fixed and known integer number, and $(\xi_t)$ a white noise. Note that given the functions $a_1,\ldots, a_T$, one may even build a periodic sequence $(a_t)_{t\in \Z}$ through the relation $a_{t+T}=a_t$.  \\
The choice of such extension of the tvAR$(1)$ processes is relative to modelling considerations: for instance, in the climatic framework, \cite{DHP} considered models of air temperatures where the function of interest writes as the product of a periodic sequence by a locally varying function. This choice provide an interesting extension of more classical periodic models of air temperature such as those proposed in \cite{LHBS}.\\
Other periodic representation for locally stationary processes can also be found in for instance in the paper \cite{pd}, but the seasonal component is treated as an additive deterministic trend and is not included in the dynamic of the process, which is the case for model \eqref{AR1}.\\
~\\
We then study non-parametric estimators $\widehat a_s(u)$, for $s=1,2,\ldots,T$, $u\in(0,1)$  from an observed trajectory $(X^{(n)}_1,\ldots, X^{(n)}_{nT})$. We consider kernel-based estimators which are naturally induced from covariance relationships satisfied by the process (see Section \ref{estim}). Central limit theorems are established for these estimators under some regularity conditions on the functions $a_s(\cdot)$ for $s=1,2,\ldots, T\,$. The results are only obtained by assuming second-order moments on the white noise $(\xi_t)$. This is a main improvement with respect to usual limit theorems on locally-stationary processes which are obtained with the assumption that any moment exists for $(\xi_t)$. This is due to the new ideas developed in our proof which combines a central limit theorem for martingale increment arrays as well as an embedding in an Orlicz space (see details in Section \ref{proofs}).       \\
The obtained convergence rate is optimal with respect to the minimax rate up to a logarithmic term. Simulations based on Monte-Carlo experiments illustrate the accuracy of the estimators. An application to real-life data, {\it i.e.} monthly average temperature readings in London from 1659 to 1998, shows the interest of using our new model \eqref{AR1}.

~\\
This paper is also a first step concerning new results for new class of non-stationary processes. 
Indeed, we can extend the definition \eqref{AR1} to processes $(X_t^{(n)})$ such as:
\begin{eqnarray}  
\label{chaineinfloc} 
X^{(n)}_t&=&F_t\left(\frac tn,\xi_t,Z_t;X^{(n)}_{t-1},X^{(n)}_{t-2},\ldots\right),\qquad 1\le t\le n,
\end{eqnarray}
where $(Z_t)$ is a sequence of i.i.d. random vectors modelling for instance exogenous inputs.  This more tough case is deferred to forthcoming papers.

Other time-varying models with an infinite memory may also be processed as GARCH-type models (see for instance \cite{dgs}). 
Remark also that \cite{ft} introduced INGARCH-models. Those models are GLM models; non-stationary versions of which also may be considered. They will be considered in further works.\\
~\\
The structure of the paper is as follows. In Section \ref{estim}, we define and study asymptotic properties of non-parametric estimators for the process \eqref{AR1}. Section \ref{simu} provides the results of some Monte-Carlo experiments and real-life data application, while the proofs are reported in Section \ref{proofs}.

\section{Asymptotic normality of a non-parametric estimator for periodic tvAR(1) processes}\label{estim}
\subsection{Definition and first properties of the process}
Denote classically $\N=\{0,1,\ldots \}$ and $\N^*=\{1,2,\ldots \}$. Here we consider  $T\in\N^*$ a fixed and known period. We will write $s\equiv t [T]$ if $t-s$ is a multiple of $T$.\\
The paper is dedicated to the  simplest case $X=(X^{(n)}_t)_{1\leq t \leq nT, \, n \in \N}$, of a $T-$periodic locally stationary  AR$(1)-$process, defined in \eqref{AR1}
where $X^{(n)}_0=X_0$ with $\E(X^2_0)<\infty$. Here  $(\xi_t)_{t\in \N}$ is a sequence of i.i.d. r.v.s satisfying $\E(\xi_t)=0$ and $\v (\xi_t)=\sigma^2$ for any $t \in \N^*$, with $(\xi_t)_t$ independent of $X_0^{(n)}$. \\
The functions $(a_s(\cdot))_{1 \leq s \leq T}$, $[0,1]\to  \R$ are supposed to satisfy some regularity. Hence, we provide the forthcoming definition usually made in a non-parametric framework:
\begin{Def}\label{Crho}
For $\rho >0$, we denote $\lceil\rho\rceil \in\N$ the largest integer  such that  $\lceil\rho\rceil<\rho$. A function $f:x\in \R \mapsto f(x) \in \R$ is said to belong to the class ${\cal C}^\rho({\cal V}_u )$ where ${\cal V}_u$ is a neighbourhood of $u\in \R$, if $f \in {\cal C}^{\lceil\rho\rceil} ({\cal V}_u )$  and if $f^{(\lceil\rho\rceil)}$ is a $\,(\rho-\lceil\rho\rceil)$-H\"olderian function, {\it i.e.} there exists $C_f \geq 0$ such as 
$$
\big | f^{(\lceil\rho\rceil)}(u_1) -f^{(\lceil\rho\rceil)}(u_2) \big | \leq C_f \, |u_1-u_2|^{\rho-\lceil\rho\rceil},\quad \mbox{for any}~u_1,u_2 \in {\cal V}_u.
$$
\end{Def}
In case $\rho$ is an integer we simply assume that  $f^{(\rho)}$ exists and is a continuous and bounded function  on the neighbourhood of $u$. As a consequence we specify the assumptions on functions $(a_t)$ using a fixed positive real number $\rho>0$: \\
~\\
{\bf Assumption (A$(\rho)$):} The functions $\{a_t(\cdot);\,{t\in \N}\}$ are such as:
\begin{enumerate}
\item (Periodicity) There exists $T \in \N^*$ such that $ a_t(v)= a_{t+T}(v)$ for any $(t,v) \in \Z\times [0,1]$. 
\item (Contractivity) There exists $\alpha=\sup_{\{t\in \Z,\, v \in [0,1]\}} |a_t(v)|<1$. 
\item (Regularity)  For any $t \in \Z$, assume that $a_t \in {\cal C}^\rho$.
\end{enumerate}
\begin{Rem}
Quote that $T=1$ corresponds to a non-periodic case and $(X^{(n)}_t)$ is then a usual tvAR(1) process defined in \eqref{tvARMA}.
\end{Rem}
First it is clear that the conditions on functions $(a_s)$ ensure the existence of a causal linear process $(X_t^{(n)})_{1\leq t \leq nT}$ for any $n \in \N$ satisfying \eqref{AR1}.  More precisely, we obtain the following moment relationships:
\begin{prop}\label{lem1}
Let $X=(X^{(n)}_t)_{1\leq t \leq nT,\, n \in \N^*}$ satisfy \eqref{AR1} under Assumption (A$(\rho)$) with $\rho\ge 1$. Then for some convenient constant $c>0$,
\begin{enumerate}
\item For any $n \in \N^*$ and $1\leq t \leq nT$, $\big |\E \big (X^{(n)}_t \big )\big | \leq \alpha^t \, \big | \E (X_0) \big |$.
\item Let $s \in \{1,\ldots,T\}$. There exists functions $\gamma^{(2)}_s\in{\cal C}^{\rho}([0,1])$ such as if  $t \in \{[c\log n],\ldots,nT\}$ and  $t\equiv s~[T ]$:
\begin{multline}\label{vt2} 
\hspace{-0cm}\E \big ((X^{(n)}_t )^2 \big  ) =\gamma_s^{(2)}(\frac t{nT})+ {\cal O}\big ( \frac 1 n \big ), \\ 
\mbox{with} \quad \left \{ 
\begin{array}{ccl}\displaystyle 
  \gamma_s^{(2)}(v)&=&\displaystyle \sigma^2  \, \frac{ 1+ \sum _{i=0}^{T-2} \beta_{s,i}(v)}{1-\beta_{s,T-1}(v) }, \\
  \displaystyle 
 \beta_{s,i}(v)&=&\prod _{j=0}^{i} a^2_{s-j}(v)  \leq \alpha^{2i}<1 .
\end{array} \right . 
\end{multline}
\item Assume $\E (\xi_0^4)=\mu_4<\infty$ and $\E (\xi_0^3)=0$ (this holds e.g. if $\xi_0$ admits a symmetric distribution).\\ For $s \in \{1,\ldots,T\}$, there exist functions $\gamma^{(4)}_s\in{\cal C}^{\rho}([0,1])$ such as, for   $t \in \{[c\log n],\ldots,nT\}$  with  $t\equiv s~[T ]$,
\begin{multline}\label{vt4} 
\E \big ((X^{(n)}_t )^4 \big  ) =\gamma_s^{(4)}(\frac t{nT})+ {\cal O}\big ( \frac 1 n \big ),\\
\quad \mbox{with} \quad  \left \{ 
\begin{array}{ccl}
  \gamma_s^{(4)}(v)&=& \displaystyle \big ( \mu_4+ 6\sigma^2\gamma_s^{(2)}(v)-6\sigma^4 \big )  \frac{1+ \sum _{i=0}^{T-2} \delta_{s,i}(v)}{1-\delta_{s,T-1}(v) }, \\
 \delta_{s,i}(v)&=&\prod _{j=0}^{i} a^4_{s-j}(v) \leq \alpha^{4i}<1. 
\end{array} \right . 
\end{multline}

Moreover, 
  for any $(t,t')\in \{[c\log n],\ldots,nT\}^2$ with $t>t'$,
\begin{equation} \label{cov2}
\cov \big ((X^{(n)}_t )^2,(X^{(n)}_{t'} )^2 \big )=  \Big ( \gamma_{s'}^{(4)}\big (\frac {t'}{nT} \big ) + {\cal O}\big ( \frac 1 n \big ) \Big) \prod_{i=1}^{t-t'} \,a^2_{t'+i}(\frac {t'+i} n).
\end{equation}
\end{enumerate}
\end{prop}
We will now assume $X_0=0$. \\
In addition of the previous proposition, another relation can be easily established. Indeed, for $t\in \{1,2,\ldots,nT\}$, with $s=t~[T]$, by multiplying \eqref{AR1} by $X_{t-1}^{(n)}$ and taking the expectation:
\begin{eqnarray}\label{a}
a_{t}\Big(\frac t{nT}\Big)=a_{s}\Big(\frac t{nT} \Big)  =\frac{\E \big (X_tX_{t-1}\big)} {\E \big (X_{t-1}^2\big)}.
\end{eqnarray}
The relation \eqref{a} is at the origin of the definition of the following non-parametric estimators of the functions $a_s(\cdot)$.  
\subsection{Asymptotic normality of the estimator}
Assume that the sample $(X_1,\ldots,X_{nT})$ is observed for some $n\ge1$; this condition entails a reasonable loss of at most $T$ data and allows us for a more comprehensive study.\\ For each $ s\in \{1,\ldots,T\}$, we define $
I_{n,s}=\big \{s,s+T,\ldots, s+(n-1)T \big \}
$, a set  with  $\# I_{n,s} =n$. 
Now \eqref{a} writes:
$$
a_{s}\Big(\frac t{nT}\Big)  =\frac{\E \big (X_tX_{t-1}\big)} {\E \big (X_{t-1}^2\big)}, \qquad \forall  t\in I_{n,s} .
$$
A convolution kernel $K:\R \to \R$  will be required in the sequel and it satisfies one of both the following assumptions:\\
~\\
\noindent{\bf Assumption $(K)$:} Let $K:\R \to \R$ be a Borel  bounded  function such that:
\begin{itemize}
\item $\displaystyle \int_\R K(t)dt=1$ and $K(-x)=K(x)$ for any $x \in \R$;
\item there exists $\beta>0$ such as $\lim_{|t|\to + \infty} e^{\beta \, |t|} K(t)=0$.  \\
\end{itemize}
\noindent{\bf Assumption $(\widetilde K)$:} Let $K:\R \to \R$ be a Borel  bounded function such that:
\begin{itemize}
\item $\displaystyle \int_\R K(t)dt=1$ and $K(-x)=K(x)$ for any $x \in \R$;
\item there exists some $B>0$ such as $K(t)=0$, if $|t|>B$.  \\
\end{itemize}
Typical examples of kernel functions are $K(t)=(2 \pi)^{-1/2} e^{-t^2/2 }$ and $K(t)=\frac 1 2 \, \1 _{ [-1,1 ]}(t)$
satisfying respectively Assumptions $(K)$ and $(\widetilde K)$. Note also the $K\ge0$ would exclude dealing with a regularity $\rho>2$.\\
~\\
For $r\geq 1$, we also specify another condition satisfied by such a function: \\
~\\
\noindent{\bf Assumption} ker$(r)$: Let $K:\R \to \R$ be a Borel bounded function such that:
\begin{itemize}
\item $\int_\R (|x|^{r}+1) \,|K(x)|\, dx<\infty$ and $\int_\R x^{p} K(x)\, dx=0$, if $p\in \{1,2, \ldots, \lceil r \rceil -1 \}$;
\item $\|K \|_\infty=\sup_{x \in \R} |K(x)| <\infty$ and $\Lip (K)=\sup_{x \neq y} \frac {|K(x)-K(y)|}{|x-y|}<\infty$.  \\
\end{itemize}
Assume that a sequence of positive bandwidths  $(b_n)_{n\in \N}$ is chosen in such a way that 
$$
\lim_{n\to\infty}b_n=0,\qquad \lim_{n\to\infty}nb_n=\infty.
$$
Now, keeping in mind the expression \eqref{a} and following the same ideas as with Nararaya-Watson estimator (see \cite{nad} and \cite{wats}), for $s \in \{1,\ldots,T\}$ and $u\in(0,1)$, we set
\begin{eqnarray}\label{est}
\widehat a^{(n)}_s(u)=\frac{\widehat N^{(n)}_s(u)}{\widehat D^{(n)}_s(u)} ,~ \mbox{with}  \left \{ \begin{array}{ccl}\displaystyle 
\widehat N^{(n)}_s(u)&\hspace{-2mm}=&\hspace{-2mm} \frac1{nb_n}\sum_{j\in I_{n,s}}
K\Big(\frac{\frac j{nT}-u} {b_n}\Big)X_jX_{j-1} ,\\
\displaystyle 
\widehat D^{(n)}_s(u)&\hspace{-2mm}=& \hspace{-2mm}  \frac1{nb_n }\sum_{j\in I_{n,s}}
K\Big(\frac{\frac j{nT}-u} {b_n}\Big) X_{j-1}^2 .
\end{array} 
\right . 
\end{eqnarray}
Since extremities are omitted we avoid the corresponding edge effects due to the fact that at the extremities, summations are not considered over a symmetric interval of times containing $nu$. The case $u=0$ does not make any contribution while the case $u=1$
corresponds with simple periodic behaviours and such results should be found in \cite{LHBS}. \\
~\\
Using essentially a martingale central limit theorem (the steps of the proofs  are precisely detailed   in Section \ref{proofs}), we obtain:
\begin{theo}\label{theo1}
Let $0<\rho\leq 2$ and Assumption (A$(\rho)$), let $K$ satisfy Assumption $(K)$ or $(\widetilde K)$ as well as Assumption ker$(\rho\vee 1)$. Then, for a sequence $(b_n)_{n\in \N}$ of positive real numbers such as $\lim_{n\to\infty}b_n \, n^{\frac1{1+2(\rho\wedge 1)}} = 0$,
\begin{equation}
\sqrt{nb_n} \big  ( \widehat a_s(u)-a_s(u) \big ) \limiteloin   {\cal N}\Big ( 0 \,  , \frac{ \sigma^2}{ \gamma^{(2)}_s(u)}\int_ \R K^2(x) \, dx \Big ),\label{TLC1} 
\end{equation}
for any $u \in (0,1)$, $s\in \{1,\ldots,T\},$ with $\displaystyle \gamma^{(2)}_s(u)=\sigma^2\, \frac{1+ \sum _{i=0}^{T-2} \beta_{s,i}(u)}{1-\beta_{s,T-1}(u) }.$ 
\end{theo}
Note that for $\rho\leq 1$ the classical optimal semi-parametric minimax rate is reached. 

This is not the case if $\rho \in (1,2]$. In that case, another moment condition is  needed  in order to improve the convergence rate of $\widehat a_s(u)$.
\begin{theo}\label{theo2}
Let $1\leq \rho\leq 2$ and Assumption (A$(\rho)$), let $K$ satisfy Assumption $(K)$ or $(\widetilde K)$ as well as Assumption ker$(\rho)$. Moreover, suppose that $\E  |\xi_0|^\beta <\infty$ with $ \displaystyle \beta= 4-\frac{2\rho}{5\rho-4} \in \Big[2,\frac {10}3\Big]$ (Note that $\beta=2$ if $\rho=2$) and that $\xi_0$ admits a symmetric distribution.
 Then \eqref{TLC1} holds for a sequence $(b_n)_{n\in \N}$ of positive real numbers such as $b_n \, n^{\frac 1{2\rho+ 1}}  \limiten 0$.
 \\
 Moreover in case $\rho=2$ and if $ b_n=c \,n^{-\frac15}$ then the central limit still holds but the limit distribution is now  non-centred: $${\cal N}\Big(\mu(u) \, ,  \, \frac{\sigma^2}{\gamma^{(2)}_s(u)} \int_ \R K^2(x) \, dx \Big )$$ with $\displaystyle \mu(u)=\frac{c^{\frac52}}{\gamma^{(2)}_s(u)}\Big(\frac12a_s''(u)\gamma_s^{(2)}(u)+ a_s'(u)(\gamma_s^{(2)})'(u)\Big)\int_\R z^2 K(z)\,dz $.
\end{theo}
\begin{Rem}
Optimal window widths write as  $b_n \sim cn^{-\frac 1{2\rho+1}} $ thus the above result holds with a suboptimal window width.
Moreover the symmetry assumption is discussed in Remark \ref{mom3}.
Now for the case $\rho=2$ in case the derivatives of $a_s$ are regular around the point $u$, then the optimal window width actually may be used and the central limit theorem again holds with a non-centred  Gaussian limit.
\\
Quote that the proposed normalisation yields the standard minimax rates $n^{-\frac\rho{2\rho+1}}$, in the case of compactly supported symmetric kernel (a $(\log n)-$loss is observed for the Gaussian kernel); the obtained rates are in probability and further work is needed to prove that this is the minimax $\L^2-$rate. 
\\
Moreover for large $T$  the convergence rate is degraded with a factor $T^{\frac\rho{2\rho+1}}$ since the sample size is $N=nT$ and thus $n=N/T$.
\end{Rem}
\begin{Rem}
Of course, if $T=1$, Theorems \ref{theo1} and \ref{theo2} hold, which provide another minimax estimation of the function $u  \mapsto a_1(u)$ ($u \in [0,1]$) requiring sharper moment and regularity conditions than the ones proposed in Theorem 4.1 of \cite{dahl2009}. 
\end{Rem}
\begin{Rem}
If $T$ is unknown we better consider an $N$-sample and set  $n=[N/T]$, the proof of previous central limit theorem \ref{theo1} provides an approach for estimating this period $T$. First fix $T_{\max} \geq 2$ (typically $T_{\max}=12$ for monthly data). Then, for each $1\leq \tau \le T_{\max}$, we define an estimator $\widehat a_s^{(\tau)}(u)$ for any $1\leq s \le \tau$ and $u \in (0,1)$. 
It is clear that when $\tau$ is not a multiple of $T$, then the sums in \eqref{est} that are done on the set $I_{n,s}$, which depends on $\tau$, is now a sum involving other $a_k$ with $k\ne s$. As a consequence, $\widehat a_s(u)$ is not a convergent estimator of $a_s(u)$. \\
Then, using a classical cross-validation, for each $1\le \tau \le T_{\max}$, we compute
$$
\widehat {CV}(\tau)=\sum_{j=2}^{N}\Big (X^{(n)}_j-\widehat a_j^{(\tau)}\big (\frac j N\big )X^{(n)}_{j-1} \Big )^2.
$$
Finally, define $\widehat T$ as the smallest value such as
\end{Rem}
\vspace{-3mm}
$$
\widehat T=\text{Arg}\! \!\! \! \! \! \! \min_{1 \le \tau \le T_{\max}} \widehat{CV}(\tau).
$$
\begin{Rem}
The central limit theorem \ref{theo1} naturally provides a test statistics $\widehat A_s$ for solving the test problem: $H_0:$ $a_s(u)=c_a$ versus $H_0:$ $a_s(u)\neq c_a$, where $c_a \in (0,1)$. Indeed, from \eqref{TLC1}  and Slutsky Lemma we deduce:
$$
\sqrt{nb_n \,\int_ \R K^2(x) \, dx } \, \sqrt {\frac{1+ \sum _{i=0}^{T-2} \prod _{j=0}^{i} \widehat a^2_{s-j}(u)}{1-\prod _{j=0}^{T-1} \widehat a^2_{s-j}(u)}}\Big ( \widehat a_s(u)-a_s(u) \Big ) \limiteloin {\cal N}\big ( 0 \,  , 1\big ).
$$
Then if we consider $$
\widehat A_s=\sqrt{nb_n \,\int_ \R K^2(x) \, dx } \, \sqrt {\frac{1+ \sum _{i=0}^{T-2} \prod _{j=0}^{i} \widehat a^2_{s-j}(u)}{1-\prod _{j=0}^{T-1} \widehat a^2_{s-j}(u)}}\Big ( \widehat a_s(u)-c_a \Big ),
$$
 this provides a natural  statistics test with usual standard Gaussian quantile as asymptotic threshold.
\end{Rem}
\section{Monte-Carlo experiments and an application to climatic data} \label{simu}
\subsection{Monte-Carlo experiments}
In this section, numerous Monte-Carlo experiments have been made for studying the accuracy of the new non-parametric estimator $\widehat a_s( \cdot)$. 

\begin{figure}
\begin{center}
\caption{\label{Fig1}\textit{Graph of the function $a_1^{(2)}$ and an example of its estimation (for $n=1000$).}}
\includegraphics[width=10cm,height=5cm]{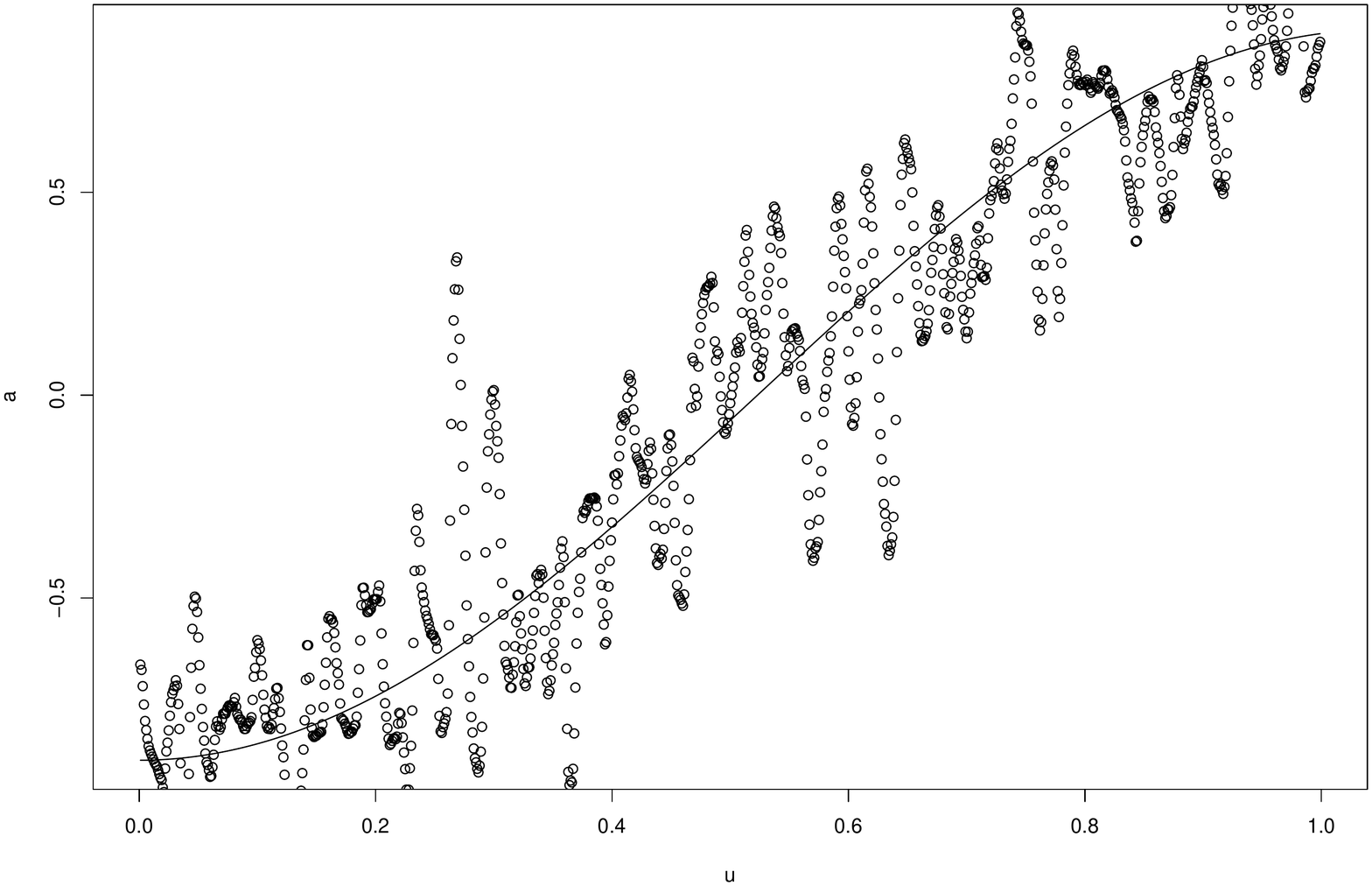} 

\end{center}
\end{figure}

\begin{figure}
\begin{center}
\caption{\label{Fig2}\textit{Graph of the chosen function $a_1^{(0.8)}$.}}
\includegraphics[width=10cm,height=5cm]{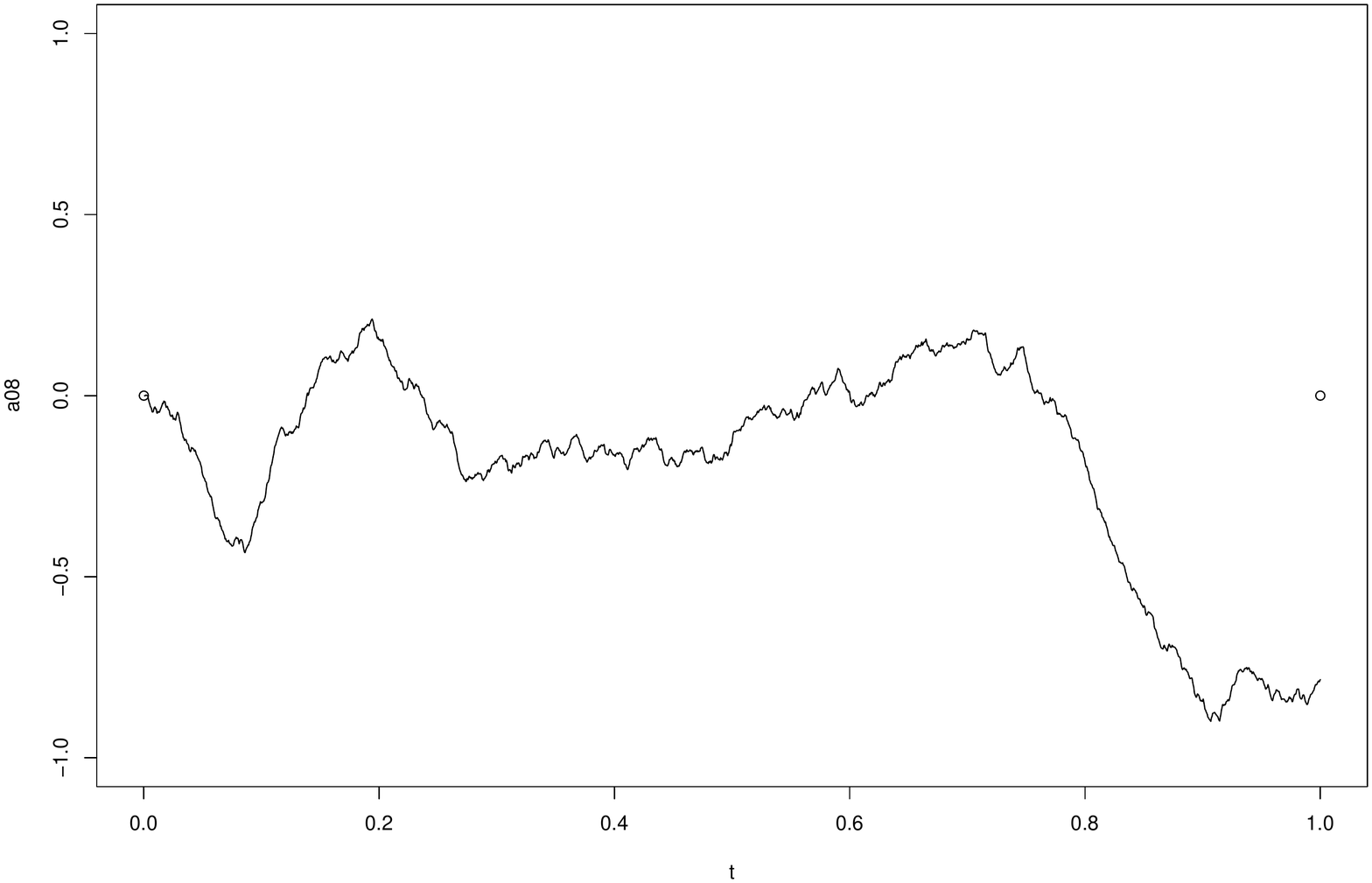}

\end{center}
\end{figure}

Firstly, we considered $3$ typical functions $[0,1] \to[-1,1]$, $ a^{(\rho)}_s(u) \in {\cal C}^\rho([0,1])$ and such as $\sup_{u\in[0,1],s\in \N} |a^{(\rho)}_s(u)|\leq \alpha <1$:
\begin{itemize}
\item For $\rho=2$, we choose $\displaystyle a_s^{(2)}(u)=0.9 \, \cos\big(2\pi\frac{ ns}T\big )\cos(3u)$. Figure \ref{Fig1} exhibits the graph of the function $a_1^{(2)}$ and an example of its estimation (for $n=1000$);
\item  For $\rho=1.5$, we choose $\displaystyle a_s^{(1.5)}(u)=0.9  \,\cos\big(2\pi\frac{ ns}T\big )\frac {\int_0^uW_t(\omega) \,dt} {\sup_{x\in[0,1]}|W_x(\omega)|} $ where $(W_t)_{t\in[0,1]}$ is an observed trajectory of a Wiener Brownian motion;
\item  For $\rho=0.8$, we choose $\displaystyle a_s^{(0.8)}(u)=0.9 \,\cos\big(2\pi\frac{ ns}T\big ) \frac {B_{0.8}(\omega,u) } {\sup_{x\in[0,1]}|B_{0.8}(\omega,x)|} $ where $B_{H}(\omega,t))_{t\in[0,1]}$ is an observed trajectory of a fractional  Brownian motion with Hurst exponent $H=0.8$ (Figure \ref{Fig2} exhibits the graph of this chosen function $a_1^{(0.8)}$). It is well known that a trajectory of a fractional  Brownian motion with Hurst exponent $H \in (0,1)$ is almost surely $\alpha$-H\"oderian for any $\alpha<H$;
\item  For $\rho=0.5$, we choose $\displaystyle a_s^{(0.5)}(u)=0.9 \,\cos\big(2\pi\frac{ ns}T\big ) \frac {W_u(\omega) } {\sup_{x\in[0,1]}|W_x(\omega)|} $ where $(W_t(\omega))_{t\in[0,1]}$ is an observed trajectory of a Wiener Brownian motion. 
\end{itemize} 
We also consider two ``typical'' kernels:
\begin{itemize}
\item A bounded supported kernel, the well-known Epanechnikov kernel defined by $K_E(x)=\frac 3 4 \, (1-x^2)\, \1 _{\{|x|\leq 1\}}$, which is known to minimize the asymptotic MISE in the kernel density estimation frame;
\item The unbounded supported Gaussian kernel with $K_G(x)=\frac 1 {\sqrt{2\pi}}  \exp\big (-\frac {x^2}2 \big )$.
\end{itemize}
We considered the cases $n=100,\, 200, \, 500$ and $1000$, and we fixed $T=2$. Finally $1000$ independent replications of $(X^{(n)})$ are generated with two different cases of innovations $(\xi_t)$:
\begin{itemize}
\item Firstly, the  case where the probability distribution of $\xi_0$ is a Gaussian ${\cal N}(0,4)$ distribution, then $\E |\xi_0|^4<\infty$ and therefore Theorem \ref{theo1} holds for $\rho=0.5$ and Theorem \ref{theo2} holds for $\rho=1.5$ and $\rho=2$.
\item Secondly, the  case where the probability distribution of $\xi_0$ is a Student $t(3)$ (with $3$ degrees of freedom) distribution implying $\E |\xi_0|^\beta<\infty$ for any $\beta<3$ but $\E |\xi_0|^3=\infty$. Then if $\rho=0.5$,  Theorem \ref{theo1} holds  but if $\rho=1.5$ and $\rho=2$, Theorem \ref{theo2} does not hold.
\end{itemize}
Finally, for each $n$, each functions $a_s^{(\rho)}$ and kernel $K$, and each probability distributions of $\xi_0$, we present the results computed from $1000$ replications and the following methodology:
\begin{enumerate}
\item For each replication $j$, we defined $b_n=n^{-\lambda}$ with $\lambda=0.10,\, 0.11,\ldots,0.80$, $(u_i)_{1\leq i \leq 99}=0.01, \,0.02,\ldots,0.99$, $s=1,2,\ldots,T$, and the estimators $\widehat a_s(u_i)$ are computed.
\item For each replication $j$ and each $\lambda=0.10,0.11,\ldots,0.80$, an estimator of the $MISE$ is computed:
$$\widehat {MISE}_s(\lambda)=\frac 1 {99} \, \sum_{i=1} ^{99} \big (\widehat a_s(u_i)-a_s(u_i) \big )^2.
$$
\item For each replication $j$, we minimised an estimator of the global square root of MISE:
$$
\widehat \lambda_j=\mbox{Arg}\!\!\!\! \!\! \! \min_{0.1\leq \lambda \leq 0.8}   \sum_{s=1}^T  \sqrt{\widehat {MISE}_s(\lambda)}
$$
\item Then we computed $\overline \lambda=\frac 1 {1000} \, \sum_{j=1}^{1000} \widehat \lambda_j$ over all the replications.
\item Finally, we computed the estimator of the minimal global square root of MISE, 
$$
\overline{MISE}^{1/2}=\frac 1 {1000} \,\sum_{j=1}^{1000}\sum_{s=1}^T  \sqrt{\widehat {MISE}_s(\widehat \lambda_j)}.
$$ 
\end{enumerate}
As a consequence, $\overline \lambda$ and $\overline{MISE}^{1/2} $ are two interesting estimators relative to Theorems \ref{theo1} and  \ref{theo2}. The first one specifies the link between the choice of an optimal bandwidth $b_n$ qnd the regularity $\rho$ of the functions $a_s(\cdot)$. The second one measures the optimal convergence rate of the estimators $\widehat a_s(\cdot)$ to $a_s(\cdot)$. All the results are printed in Tables \ref{Table1} and \ref{Table2}. \\
Moreover, for exhibiting the asymptotic normality of the estimators provided in the central limit theorem \eqref{TLC1}, we draw in Figure the histograms of $\widehat a^{(2)}_s(u)$ for $u=0.25,\, 0.5$ and $0.75$ from $10000$ independent replications for $n=5000$. We  also used a Jarque-Bera test to confirm the Gaussian asymptotic distribution since the p-values of this test are successively: $p-value=0.105$, $0.927$ and $0.345$. Hence, the asymptotic normality of the estimator seems to be attested by Monte-Carlo experiments.
~\\ 
\begin{figure}
\begin{center}
\caption{\label{Fig3}\textit{Histograms of $\widehat a^{(2)}_s(u)$ for $u=0.25,\, 0.5$ and $0.75$ from $10000$ independent replications}}
\includegraphics[width=3.7cm,height=3cm]{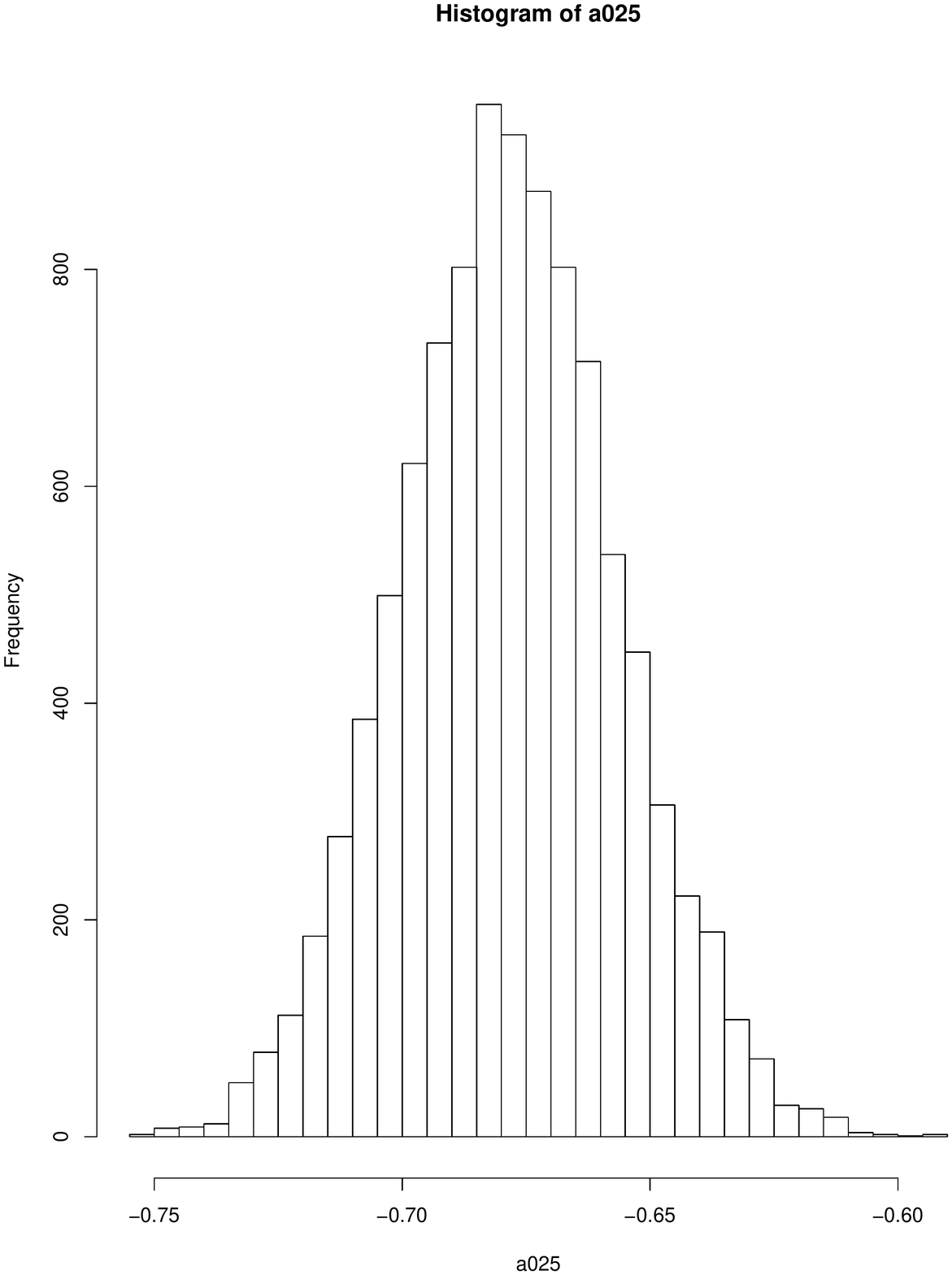}
\includegraphics[width=3.7cm,height=3cm]{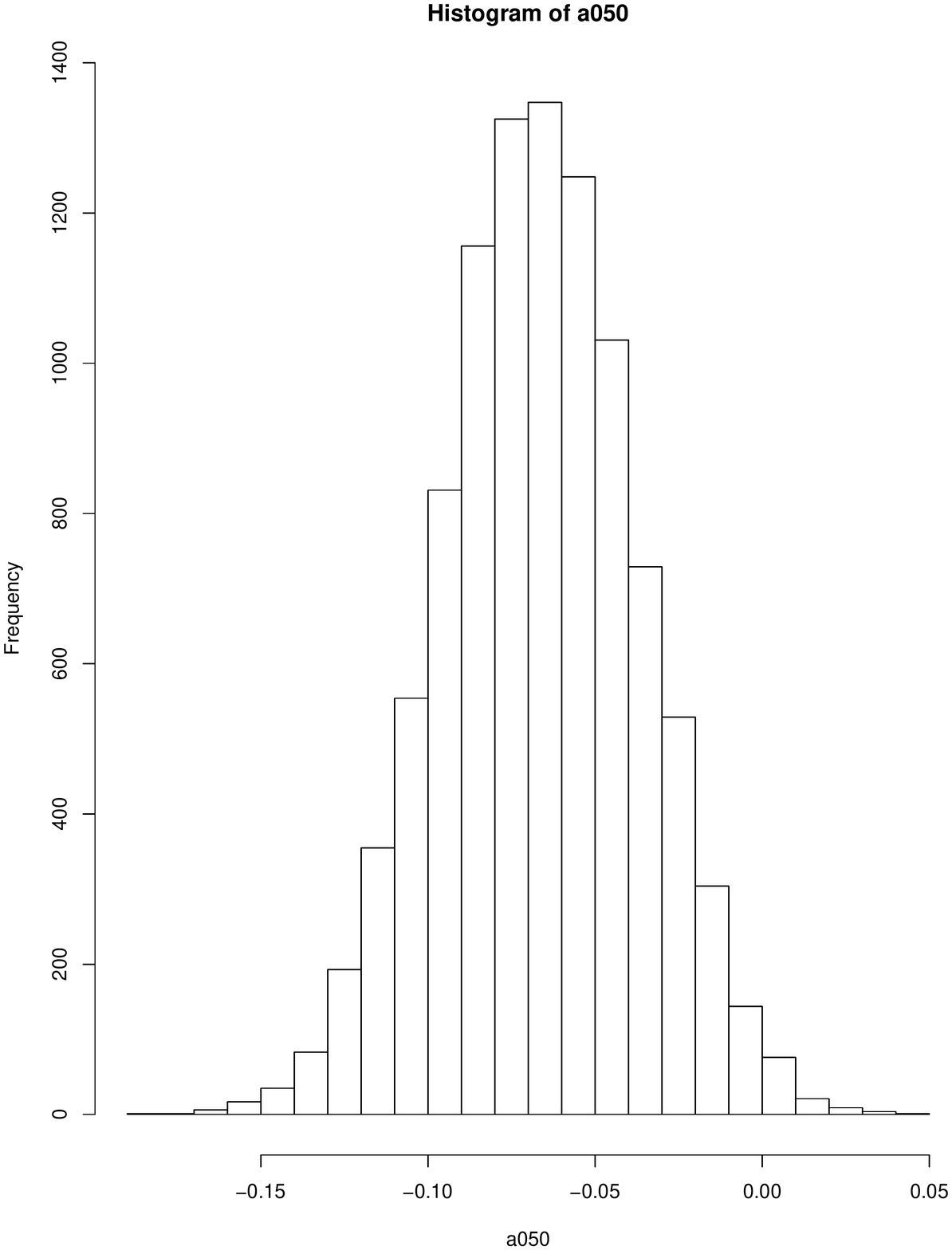}
\includegraphics[width=3.7cm,height=3cm]{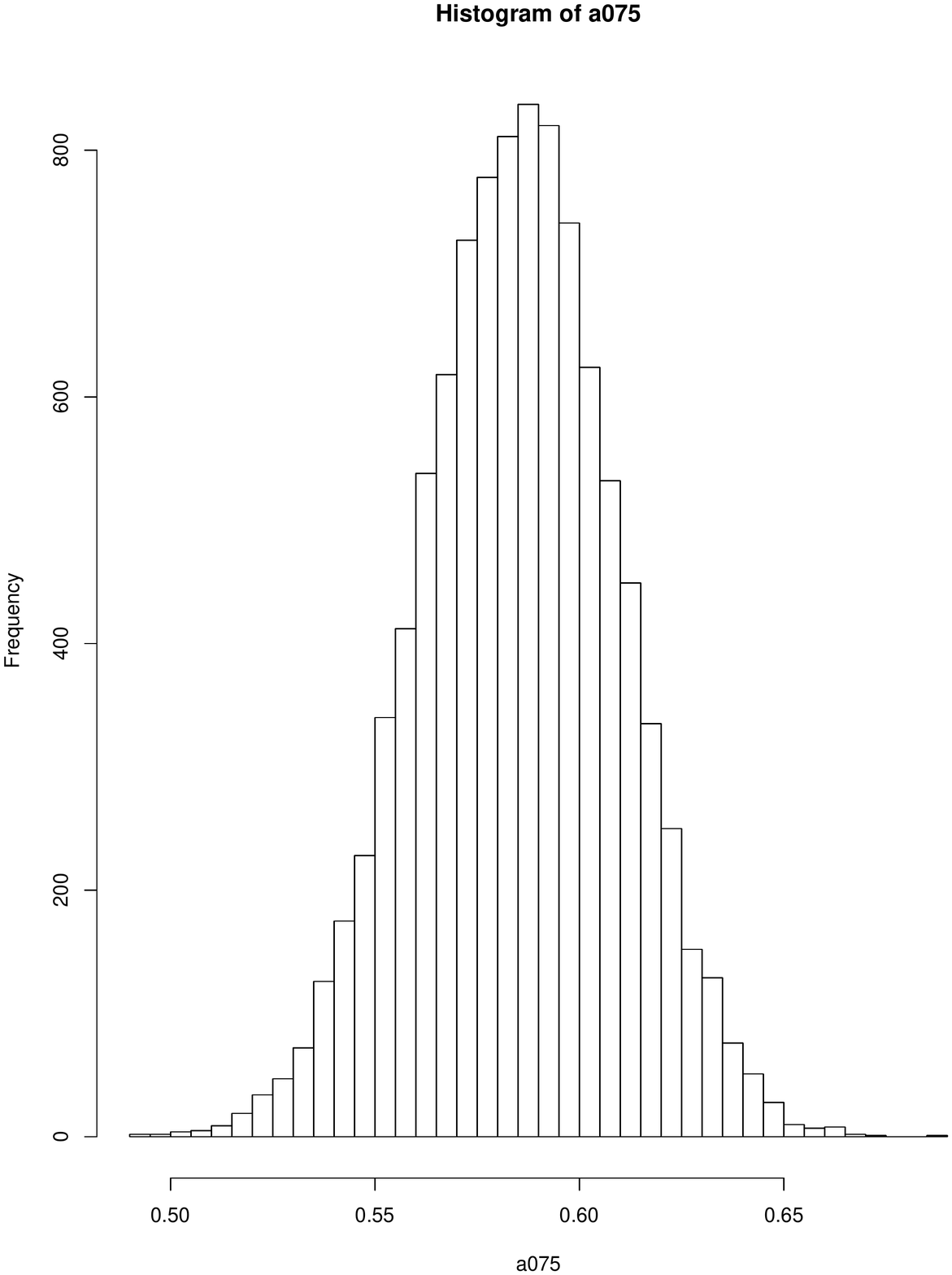}
\end{center}
\end{figure}
\begin{table}
\caption{ \label{Table1} Results of the Monte Carlo experiments providing  the accuracy of $\widehat a_s$ for the three chosen functions  the three chosen functions with $\xi_0$ following a ${\cal N}(0,4)$ distribution, $1000$ independent replications are generated.}
\begin{center}
\begin{tabular}{|c|c||c|c||c|c||c|c||c|c||}
\hline
 & $a^{(\rho)}$& \multicolumn{2}{c|}{$a_s^{(2)}$} & \multicolumn{2}{c|}{$a_s^{(1.5)}$}& \multicolumn{2}{c|}{$a_s^{(0.8)}$} & \multicolumn{2}{c|}{$a_s^{(0.5)}$} \\
 \hline
 & Kernel & $K_E$ &$K_G$  & $K_E$ &$K_G$ & $K_E$ &$K_G$ & $K_E$ &$K_G$ \\
\hline\hline 
$n=100$& $\overline \lambda$ &   0.243 & 0.407& 0.283&0.450 & 0.172& 0.322 & 0.235& 0.392  \\ \hline 
&$\overline{MISE}^{1/2}$ & 0.248 & 0.239 & 0.286 & 0.282 &0.230& 0.234  &0.354& 0.353\\  \hline 
$n=200$ & $\overline \lambda$  &0.227 & 0.363& 0.278& 0.429 &0.256 & 0.392 &0.250 & 0.386\\  \hline
&$\overline{MISE}^{1/2} $&  0.185 & 0.175 &0.219  & 0.219 &0.232&0.232 &0.308&0.303\\ \hline
$n=500$ & $\overline \lambda$  & 0.234& 0.320 &0.276 & 0.399 &0.321 &0.431 &0.287 &0.406 \\ \hline
&$\overline{MISE}^{1/2} $ & 0.129 & 0.119 & 0.154 & 0.156 &0.213& 0.210 &0.256& 0.254 \\ \hline
$n=1000$ & $\overline \lambda$  & 0.240 & 0.321 &0.270&   0.384& 0.373&0.476 & 0.328&0.438\\ \hline
&$\overline{MISE}^{1/2} $ & 0.098 & 0.093 &0.124& 0.122& 0.207&0.202 & 0.226&0.221 \\
\hline
\end{tabular}
\end{center}
\end{table}
\

\noindent {\bf Conclusions of the simulations:} Firstly, and as it should be deduced from Theorem \ref{theo1} and \ref{theo2}, we observed the larger the regularity $\rho$, the smaller $\overline \lambda$ and therefore the larger the optimal bandwidth $\overline{b_n}=n^{-\overline \lambda}$, and the faster the convergence rate of $\widehat a_s$. Secondly, even if the choice of the optimal bandwidth is significantly different following the choice of the kernel (clearly smaller with the Epanechnikov kernel), the optimal convergence rate is almost the same for both the kernel. Finally, according also with Theorem \ref{theo2}, the convergence rate is clearly slower with a heavy tail distribution ($t(3)$) than with a Gaussian distribution, and this phenomenon increases when $\rho$ increases.
\begin{table}
\caption{ \label{Table2} Results of the Monte Carlo experiments providing  the accuracy of $\widehat a_s$ for the three chosen functions with $\xi_0$ following a $t(3)$ distribution, $1000$ independent replications are generated.}
\begin{center}
\begin{tabular}{|c|c||c|c||c|c||c|c||c|c||} \hline
& $a^{(\rho)}$& \multicolumn{2}{c|}{$a_s^{(2)}$} & \multicolumn{2}{c|}{$a_s^{(1.5)}$}& \multicolumn{2}{c|}{$a_s^{(0.8)}$} & \multicolumn{2}{c|}{$a_s^{(0.5)}$} \\
 \hline
 & Kernel & $K_E$ &$K_G$  & $K_E$ &$K_G$ & $K_E$ &$K_G$ & $K_E$ &$K_G$ \\
\hline\hline 
$n=100$& $\overline \lambda$ &   0.226 & 0.394& 0.267&0.430 & 0.161& 0.295  & 0.220& 0.360 \\ \hline
&$\overline{MISE}^{1/2}$ & 0.341 & 0.320 & 0.350 & 0.340 &0.311& 0.309 &0.418& 0.405  \\  \hline 
$n=200$ & $\overline \lambda$  &0.207 & 0.343& 0.259& 0.402 &0.231 & 0.355 &0.225 & 0.362\\ \hline
&$\overline{MISE}^{1/2}$&  0.261 & 0.258 &0.281  & 0.287 &0.296&0.293 &0.353&0.346\\ \hline
$n=500$ & $\overline \lambda$  & 0.194& 0.304 &0.252 & 0.373 &0.286 &0.383 &0.239 &0.360 \\ \hline
&$\overline{MISE}^{1/2}$ & 0.214 & 0.201 & 0.213 & 0.217 &0.269& 0.261 &0.302& 0.296\\ \hline
$n=1000$ & $\overline \lambda$  & 0.193 & 0.321 &0.246&   0.342& 0.346&0.450 & 0.258&0.368 \\ \hline
&$\overline{MISE}^{1/2}$ & 0.166 & 0.093 &0.172& 0.181& 0.258&0.250 & 0.262&0.275  \\
\hline
\end{tabular}
\end{center}
\end{table}
\subsection{Numerical application on climatic data}
\begin{figure}
\begin{center}
\caption{\label{Fig4}\textit{Monthly average temperature readings in London from 1659 to 1998: correlogram of residual data}}
\includegraphics[width=12cm,height=4cm]{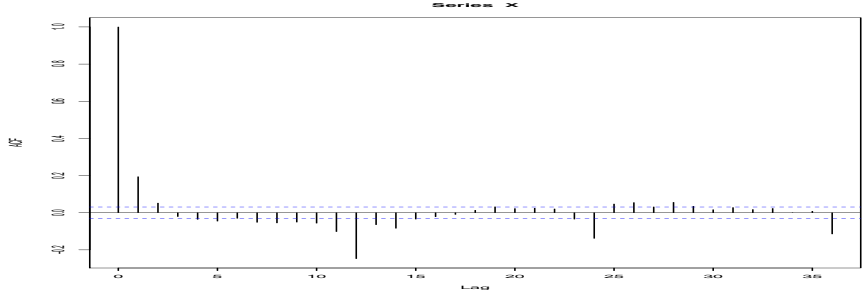}
\end{center}
\end{figure}

\begin{figure}
\begin{center}
\caption{\label{Fig5}\textit{Monthly average temperature readings in London from 1659 to 1998: $\widehat a_s(u)$ estimator for $u=0.25$ (black), $u=0.50$ (red) and $u=0.75$ (blue) in terms of $s=1,\ldots, 12$}}
\includegraphics[width=12cm,height=4cm]{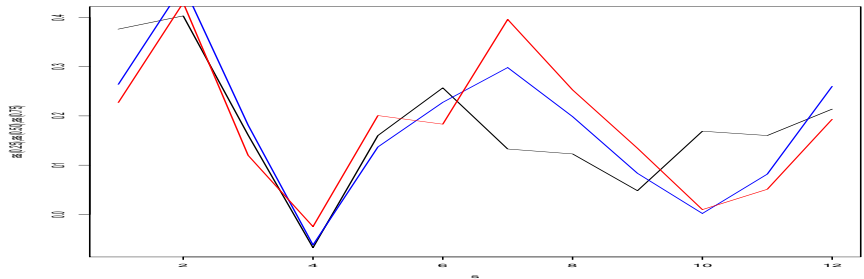}
\end{center}
\end{figure}

We also applied our model and its estimator to an example of real data, specifically the monthly average temperature readings in London from 1659 to 1998, or 340 years. Obviously in such a case one can expect that $T=12$. \\
First, we removed an additive seasonal and trend component (estimated by LOESS) from these data and considered the residual data. On these, a global correlogram (see Figure \ref{Fig4}) confirms a modelling by a process of type AR($1$) and also the presence of a periodic phenomenon of period $12$.
As a consequence we may assume that these residual data can be modelled by the model \eqref{AR1}. We then applied the $\widehat a_s(u)$ estimator for $u=0.25, \, 0.50$ and $0.75$ and $s=1,\ldots, 12$. Figure \ref{Fig5} summarizes these results and shows:
\begin{itemize}
\item The crucial interest of taking a pseudo-periodic model as we defined it in \eqref{AR1};
\item The relatively small but not negligible change in the coefficient $a_t(t)$ as a function of $t$.
\end{itemize}

\section{Proofs}\label{proofs}
We first provide the proof of Proposition \ref{lem1}. \\

{\it Proof of Proposition \ref{lem1}.} 
\begin{enumerate}
\item We have $ \E  X_1^{(n)} =a_1\Big (\frac 1 {nT} \Big ) \E (X_0 )$  and  $ \E  X_t^{(n)} =a_t\Big (\frac t {nT} \Big ) \E  X_{t-1}^{(n)})$ 
from the relation \eqref{AR1}. From Assumption (A$(\rho)$) and since $\Big  |a_1\Big(\frac 1 {nT}\Big) \Big | \leq \alpha<1$, we deduce the first item of Proposition \ref{lem1}. \\

\item Below, for ease of reading, we will omit the exponent $n$. Set  $v_t=\E \big ( X_t ^2 \big )$, and  $v=\sup_s v_{s}\in[0,+\infty]$;  also write $\alpha_t= a^2_t\big (\frac t {nT} \big )$. We have:
\begin{eqnarray}
\label{v1} v_t  = \alpha_t \, v_{t-1} + \sigma^2 
 \leq \alpha^2 v_{t-1} + \sigma^2  \leq \alpha^2 \sup_sv_s+ \sigma^2, \quad t>0
\end{eqnarray}
thus
\begin{eqnarray}
\label{bornev} \sup_s v_{s}  \leq  \frac {\sigma^2 +  v_0} {1 -\alpha}<\infty.
\end{eqnarray}

Moreover, with $\delta_t=v_t-v_{t-T}$ for any $t > T$, we have
\begin{eqnarray}
\nonumber  \delta_t & = & \alpha_t  \delta_{t-1} + (\alpha_t - \alpha_{t-T} )  v_{t-T-1},  \\
\label{bornedelta}\big | \delta_t \big |  & \leq & \alpha  \big |\delta_{t-1}\big | +  C  |\alpha_t - \alpha_{t-T} |, \quad \mbox{with}\quad  C>0,
\end{eqnarray}
from \eqref{bornev} and since for some constant $C>0$,
\begin{eqnarray}\nonumber 
\big |\alpha_t -\alpha_{t-T} \big |&=&
\Big | a^2_t \big (\frac t {nT} \big ) - a^2_{t-T} \big (\frac {t-T} {nT} \big ) \Big |
\\
&\le &2\alpha \ | a_t\big (\frac t {nT} \big ) - a_{t-T} \big (\frac {t-T} {nT} \big ) \Big |\nonumber 
\\
&\le & \frac {C} {n^{\rho\wedge1}}
\label{borneaa}
\end{eqnarray}
from  Assumption (A$(\rho)$).
As a consequence of \eqref{bornedelta}, we also obtain:
\begin{equation*}
\big | \delta_t \big | \leq \frac C {1-\alpha}\cdot \frac 1 {n^{\rho\wedge1}}+\delta_{T+1}\alpha^{t-T+1}. 
\end{equation*}
Thus for other constants $C,c>0$ we derive  
\begin{equation}\label{diffdelta}
\big | \delta_t \big | \leq  C'  \frac 1 {n^{\rho\wedge1}},\qquad \forall t\ge c\log n. 
\end{equation}
From now on, assume that $\rho\ge1$.\\
Now use again the definition  \eqref{AR1} of the model, and by iterating \eqref{v1}, we derive:
\begin{eqnarray*}
v_t&=&\sigma^2 + \alpha_t \big (\sigma ^2 + \alpha_{t-1} v_{t-2} \big ) \\
&= & \cdots \\
& =& \sigma^2 \Big (1 + \sum _{i=0}^{T-2} \alpha_{t} \cdots  \alpha_{t-i} \Big )+  \alpha_{t} \cdots  \alpha_{t-T+1} \, v_{t-T} \\
& =& \sigma^2 \Big (1 + \sum _{i=0}^{T-2} \alpha_{t} \cdots  \alpha_{t-i} \Big )+  \alpha_{t} \cdots  \alpha_{t-T+1} \, v_{t}+ {\cal O}\big ( \frac 1 n \big ) 
\end{eqnarray*} 
from \eqref{diffdelta}.\\ Hence,
\begin{eqnarray}\label{vt}
v_t&=&\sigma^2 \, \frac{ 1+ \sum _{i=0}^{T-2} \alpha_{t} \cdots  \alpha_{t-i}}{1-\alpha_{t} \cdots  \alpha_{t-T+1} } + {\cal O}\big ( \frac 1 n \big ).
\end{eqnarray} 
Now quoting that  $ \displaystyle \alpha_{t-j}=   a_{t-j}^2 \big (\frac {t-j} {nT} \big ) $ we set  $ \displaystyle \widetilde \alpha_{t-j}= a_{t-j}^2 \big (\frac t {nT}  \big) $ for $1\le j<T$, then since $\rho\ge1$ and from (\ref{vt}) we derive
\begin{eqnarray}\label{vttilde}
v_t=\sigma^2 \, \frac{ 1+ \sum _{i=0}^{T-2} \widetilde\alpha_{t} \cdots \widetilde \alpha_{t-i}}{1-\widetilde\alpha_{t} \cdots\widetilde  \alpha_{t-T+1} } + {\cal O}\big ( \frac 1 n \big )=\gamma_s^{(2)}(\frac{t}{nT})+ {\cal O}\big ( \frac 1 n \big ).
\end{eqnarray} 
The conclusion follows. 

\item
The proof mimics the case of $\E(X^2_t)$.  Denote  $q_t= a_t^4\big (\frac t {nT} \big )$, and $\mu_k=\E (\xi_0^k)$, for $k=1,2,3,4$. Then $\mu_1=0$ and
$$
w_t=\E (X^4_t)=\E(A_tX_{t-1}+\xi_t)^4=q_t \, w_{t-1}+4 \, \mu_3\, A_t  \, \E X_{t-1}+6 \,\sigma^2 \,  A_t^2 \,v_{t-1}+\mu_4.
$$
Since $\mu_3=0$, we have:
\begin{equation}\label{w}
w_t  = q_t  w_{t-1} + 6 \sigma^2 v_t + \mu_4-6\sigma^4 \leq \alpha^4 \, w_{t-1} + r(t),
\end{equation}
with $r(t)=6\, \sigma^2 \, v_t + \mu_4-6\sigma^4$ and this implies as previously $\sup_{t} w_t <\infty$. 
We also obtain for constants again denoted $C', C''>0$: 
\begin{equation}\label{diffdelta2}
\big |w_t-w_{t-T} \big | \leq \frac {C'}  n, \qquad \forall t\ge c\log n. 
\end{equation}
Finally by iterating \eqref{w}, we obtain:
\begin{eqnarray*}
w_t&=& q_{t} \cdots  q_{t-T+1} \, w_{t-T}+\Big (r(t) + \sum _{i=0}^{T-2} \big (q_{t} \cdots q_{t-i} \big )r(t-i-1)\Big ) \\
&=& q_{t} \cdots  q_{t-T+1} \, w_{t} + {\cal O}\big ( \frac 1 n \big )+\Big (r(t) + \sum _{i=0}^{T-2} \big (q_{t} \cdots q_{t-i} \big )r(t-i-1)\Big )
\end{eqnarray*} 
from \eqref{diffdelta2}. Hence, always following the previous case  
\begin{eqnarray*}\label{wt}
w_t&=&  \frac{ r(t) + \sum _{i=0}^{T-2} \big (q_{t} \cdots q_{t-i} \big )r(t-i-1)}{1-q_{t} \cdots  q_{t-T+1} } + {\cal O}\big ( \frac 1 n \big ) \\
&=& \frac{ r(t) + \sum _{i=0}^{T-2} \big (q_{t} \cdots q_{t-i} \big )r(t-i-1)}{1-q_{t} \cdots  q_{t-T+1} } + {\cal O}\big ( \frac 1 n \big ), 
\end{eqnarray*} 
for $t\ge C''\log n$,
and this implies \eqref{vt4} from using again the regularity of the functions $(a_i)_{1\le i\le T}$. \\
Finally, for any $t >t'$ such that $t,t' \in \{[c\log n],\ldots,nT\}$, since $(X_t)$ is a causal process and by iteration,
\begin{eqnarray*}
\cov (X^2_t,X^2_{t' })&=&  \alpha_t \, \cov (X^2_{t-1},X^2_{t'}) +0+ \cov (\xi ^2_{t},X^2_{t'})  \\
& = &\alpha_t \,\cov (X^2_{t-1},X^2_{t'})^{\ } \\
& = &   \big (\gamma_{s'}^{(4)}\big (\frac {t'}{nT} \big )+ {\cal O}\big ( \frac 1 n \big )\big )\prod_{i=1}^{t-t'} \alpha_{t'+i},
\end{eqnarray*}
where $s'\equiv t' ~[T]$ and $\displaystyle \Big |\prod_{i=1}^{t-t'} \alpha_{t'+i} \Big| \leq  \alpha^{2 |t-t'|}$. 
\end{enumerate}
This completes the proof.\findem

Now we establish a technical lemma, which we were not able to find in the past literature (even if variants of this result may be found) and that will be extremely useful in the sequel. For a bounded continuous function $c$ defined on $[0,1]$, and a kernel function $H$ (see details below),  an approximation of integral by appropriate Riemann sums yields (as for \cite{PrCh}'s estimator, see \cite{ros1} for further developments):
\begin{equation*}
\lim_{n\to\infty}\frac 1 {nb_n } \,\sum_{j\in I_{n,s}}H\left(\frac{\frac j{nT}-u} {b_n}\right)c\Big(\frac j{nT}\Big) =c(u),
\end{equation*}
where $u \in (0,1)$, $I_{n,s}=\big \{ s, s+T, \ldots, s+(n-1)T\}$ with $s \in \{1,\ldots,T\}$ and $T\in \N^*$. More precisely we would like to provide expansions of
\begin{equation}\label{Delta}
\Delta_n=\frac 1{ nb_n  }\sum_{j\in I_{n,s}}H\left(\frac{\frac j{nT}-u} {b_n}\right)c\Big(\frac j{nT}\Big) -c(u).
\end{equation}
\begin{lem}\label{Riemann}
Let $u\in (0,1)$, $\rho >0$, $c \in {\cal C}^\rho([0,1])$ a bounded function. Let $H$ satisfy ker$(\rho \vee 1)$.  Consider also a sequence of positive real numbers $(b_n)_n$ satisfying  $\lim_{n\to\infty}b_n=0$. Then,
there exists $C>0$ depending only on $\|H \|_\infty$, $\|c \|_\infty$ and $\Lip(H)$ such that for $n$ large enough 
\begin{equation}\label{bound}
\big | \Delta_n  \big | \leq C \,  \Big (\frac {A_n}{nb_n}  +b_n^{\rho} \Big ) , \mbox{with }  \Big \{ \begin{array} {ll}A_n=1&\hspace{-0.2cm} \mbox{under Assumption $( \widetilde K)$}, \\
A_n= \log (n)& \hspace{-0.2cm} \mbox{under Assumption $(K)$}. \end{array} 
\end{equation}
Finally, if $\rho\in\N^*$ we have:
\begin{equation}\label{equbiais}
\Delta_n= b_n^\rho\cdot \frac{c^{(\rho)}(u)}{\rho!}\, \int_\R z^\rho H(z)\,dz \,\big (1+o(1) \big )+{\cal O}\Big(\frac{A_n}{nb_n}\Big).
\end{equation}
\end{lem}

{\it Proof of Lemma \ref{Riemann}.}  In the sequel we will denote $h_n(v)=\frac1{b_n}  H\big(b_n^{-1} (v-u)\big)$ for $v \in \R$. Then $h_n$ is a Lipschitz function with $\Lip h_n=\frac 1 { b_n^2} \, \Lip H$. 
\begin{itemize}
\item
{\it First assume that the function $c\equiv1$ is a constant. }
Set $v_i=i(nT)^{-1}$ for $i \in \Z$. For $1\le s\le T$, we consider the sets 
\begin{eqnarray*}
K_{n,s}&=&\big \{ j \in \N,\, | v_{s+jT}-u| \leq A_n b_n, 1\le j\le n \big \} \\
&=& \N\cap [1,n] \bigcap \left [ (u-A_nb_n)n -\frac sT\, ,\, (u+A_nb_n)n -\frac sT\right],
\end{eqnarray*}
and $L_{n,s}=I_{n,s} \setminus K_{n,s}$. Then, for $n$ large enough,
\begin{eqnarray*}
\Delta_n & =& \frac 1 {n} \, \sum_{i\in I_{n,s}}h_n(v_i) -\int_\R h_n(v)\,dv \\
& =& \frac 1 {n}  \sum_{i\in K_{n,s}}h_n(v_i) -\int_\R h_n(v)\,dv+\frac 1 {n}  \sum_{i\in L_{n,s}}h_n(v_i)\\
& =& \sum_{j=[(u-A_nb_n)n -\frac sT ]+1}^{[(u+A_nb_n)n -\frac sT]}\int_{v_{s+jT}}^{v_{s+(j+1)T}}\big (h_n(v_{s+jT})-h_n(v)\big )dv \\
&+&\frac 1 {n} \, \sum_{j\in L_{n,s}}h_n(v_{s+jT}) -\int_{(u+A_nb_n)+\frac1n}^\infty h_n(v)\,dv-\int_{-\infty}^{(u-A_nb_n)+1/n} h_n(v)\,dv .
\end{eqnarray*}
But $|h_n(v_{s+jT})| \leq \frac C {b_n} \, \exp \Big (-\beta  \Big | \frac {j/n-u+s/nT}{b_n} \Big | \Big )$ from Assumption $(K)$ and using the usual comparison between sums and integrals for monotonic functions, we obtain:
\begin{eqnarray*}
\Big |\frac 1 {n} \, \sum_{j\in L_{n,s}}h_n(v_{s+jT}) \Big | &\leq &\frac C {b_n}\, \int_{|v-u| \geq A_nb_n}^\infty \exp \Big (- \beta \,  \frac {|v-u|}{b_n}\Big )\,dv \\
& \leq & 2 \, C\,  \exp(-\beta \, A_n).
\end{eqnarray*}
Thus
\begin{eqnarray*}
 \big |\Delta_n  \big | & \leq & \Lip (h_n)  \sum_{j=[(u-A_nb_n)n - \frac sT ]+1}^{[(u+A_nb_n)n -\frac sT ]}\int_{v_{s+jT}}^{v_{s+(j+1)T}}(v-v_{s+jT})dv  \\
&& + 2 \, C\,  \exp(-\beta \, A_n)+\int_{u+A_nb_n}^\infty\big | h_n(v) \big |\,dv+\int_{-\infty}^{u-A_nb_n} \big | h_n(v) \big |\,dv \\
 & \leq & \frac {\Lip (H) } {b_n^2} \,\frac {2A_nb_nn}{2n^2}  +  4 \, C\,  \exp(-\beta \, A_n) \\
& \leq & \Lip (H) \, \frac {A_n } {nb_n} +  4 \,C\,  \exp(-\beta \, A_n).
\end{eqnarray*}
 because since $u\in(0,1)$, the above indices remain in the index set $[-n, n]$ for $n$ large enough.\\
Then, if $A_n \geq \beta^{-1} \log n$ then  $\exp(-\beta \, A_n) \leq 1/n$ and we deduce \eqref{bound}.
\item
{\it We  now turn to the case of a non-constant function $c$.} First, if $\rho>0$, for $(u,v)\in (0,1)^2$ the Taylor-Lagrange formula implies:
$$
c(v)-c(u)=(v-u)c'(u)+\cdots+\frac{(v-u)^\ell}{\ell!} \, c^{(\ell)}\big (u+\lambda(v-u)\big ), 
$$ 
with $\ell =\lceil\rho\rceil$ and $\lambda \in (0,1)$. Since $c \in {\cal C}^\rho([0,T]) $, 
$$\big | c^{(\ell)}\big (u+\lambda(v-u)\big )-c^{(\ell)}(u) \big |\leq C_\rho \,\big |\lambda(v-u)  \big |^{\rho-\ell}\leq C_\rho \,\big |v-u \big |^{\rho-\ell}.
$$
Therefore,
\begin{equation}\label{TaylorHolder}
c(v)-c(u)=(v-u)c'(u)+\cdots+\frac{(v-u)^\ell}{\ell!} \, c^{(\ell)} (u)+ R(u,v), 
\end{equation}
with $|R(u,v)|\le C_\rho \, |u-v|^{\rho}$. 
Then for any $u\in (0,1)$, using Assumption ker$(\rho \vee 1)$ and especially the relation $\int z^pH(z)dz=0$ for $p=1,\ldots,\ell$,
\begin{eqnarray}
 \nonumber \Big |\int_\R h_n(v)\,c(v) \, dv - c(u) \int_\R h_n(v)\, dv  \Big | &=& \Big |\int_{-\infty}^\infty H(z)\big (c(u+b_nz)-c(u)\big )\,dz \Big | \\
 &=&\Big | \int_{-\infty}^\infty  H(z)R(u,u+b_nz) dz \Big |  \label{pnut0}\\
\nonumber &\leq &C_\rho \, b_n^{\rho} \,  C \, \int_{-\infty}^\infty e^{-\beta |z|} |z|^\rho dz  \\
\label{pnut} &\leq &C'\, b_n^{\rho}
\end{eqnarray}
with $C'>0$. 
Here we denote $k_n(v)=h_n(v) c(v)$ for $v \in [0,1]$. \\
Now, if $\rho \in (0,1)$, we have $$|k_n(v_1)-k_n(v_2)| \leq \|c \|_\infty \, \Lip(h_n) \, |v_1-v_2|+\frac {\|H\|_\infty} {b_n} \, C_\rho \, |v_1-v_2|^\rho,$$ and therefore using the previous results:
\begin{eqnarray*}
\big |\Delta_n  \big | &\leq &  \sum_{j=[(u-A_nb_n)n -\frac sT  ]+1}^{[(u+A_nb_n)n -\frac sT ]}\int_{v_{s+jT}}^{v_{s+(j+1)T}} C\, \Big (\Lip(h_n) \, |v-v_{s+jT}|\\
&& \hspace{1cm}+\frac {1} {b_n}  |v-v_{s+jT}|^\rho\Big )\,  dv  +  C\, \|c\|_\infty \,  \exp(-\beta  A_n)  \\
&&\hspace{2cm} +  \Big |\int_\R h_n(v)\,c(v) \, dv - c(u) \int_\R h_n(v)\, dv  \Big | \\
&\leq &   C \Big (\frac {A_n } {nb_n} + \frac {A_n } {n^\rho}+   \exp(-\beta A_n) +  b_n^{\rho}\Big ).
\end{eqnarray*}
from \eqref{pnut} and this implies \eqref{bound} since $nb_n\to \infty$ and therefore $n^{-\rho}$ is negligible with respect from $b_n^\rho$. \\
Now, if $\rho\geq 1$ and since $H$ and $c$ are bounded continuous Lipschitz functions, we obtain the inequality $$\Lip (k_n )\leq \|c\|_\infty \Lip (h_n)+ \frac 1 {b_n} \, \|H\|_\infty \Lip(c)<\infty.
$$ 
Then, using the same computations than previously (replace $h_n$ by $h_n \times c$),
\begin{eqnarray*}
|\Delta_n |& \leq & \sum_{j=[(u-A_nb_n)n -\frac sT  ]+1}^{[(u+A_nb_n)n -\frac sT ]}\int_{v_{s+jT}}^{v_{s+(j+1)T}} \big |k_n(v_{s+jT}) - k_n(v) \big | \, dv\\
& +& C\, \|c \|_\infty \exp(-\beta A_n) + \Big |\int_\R h_n(v)\,c(v) \, dv - c(u) \int_\R h_n(v)\, dv  \Big |\\
&  \leq & C \frac {A_n}{n b_n}\Big ( \|c\|_\infty   \Lip (H)+ b_n \,\Lip(c) \|H\|_\infty  \Big )+ C \|c \|_\infty e^{-\beta A_n}  + C' b_n^{\rho},
\end{eqnarray*}
from \eqref{pnut} and this completes the first item since $b_n$ is supposed to converge to $0$. The proof is now easily completed.
\item {\it Finally, in the case $\rho \in \N^*$}, we can use the previous case an a Taylor-Lagrange expansion of the function $c$, implying $\displaystyle R(u,v)=\frac{c^{(\rho)}(\theta)}{\rho !} \, \big | u -v \big |^\rho$ with $\theta =\lambda u+(1-\lambda)v$ and $\lambda\in [0,1]$.\\
 Then, using \eqref{pnut0} and with $\mu_u(z) \in [0,1]$, and $\displaystyle\zeta_n= \int_\R h_n(v)\,c(v) \, dv - c(u) \int_\R h_n(v)\, dv$
\begin{eqnarray*}
 \nonumber \zeta_n &=&\frac {b_n^\rho}{\rho!} \,  \int_{-\infty}^\infty  H(z)z^{\rho} c^{(\rho)}\big (u + \mu_u(z)b_n z\big )  \,dz \\
&=&\frac {b_n^\rho}{\rho!} \, c^{(\rho)}(u)  \,  \int_{-\infty}^\infty  H(z)z^{\rho} \,dz  \, \big (1+o(1) \big ) 
\end{eqnarray*}
from  Lebesgue theorem on dominated convergence. \findem
\end{itemize}

In the sequel we will denote the $\sigma$-algebra
\begin{equation}\label{F}
{\cal F}^{(s)}_{t}=\sigma\big ((\xi_{i})_{i\leq s+(t-1)T} \big ).
\end{equation}
\begin{lem}\label{mixingale}
Let $H$ satisfy Assumption ker$(1)$ and $(X^{(n)}_t)$  be a solution of \eqref{AR1} under Assumption (A$(\rho)$) with $\rho>0$. Then for any $u \in (0,1)$, and $s\in \{1,\ldots,T\}$,
$$
 \frac 1 {n b_n} \, \sum_{j=1}^n H\Big(\frac{\frac {s+(j-1)T}{nT}-u} {b_n}\Big) \big (X^{(n)}_{s+(j-1)T-1}\big )^2
 \limiteproban \sigma^2\, \frac{ 1+ \sum _{i=0}^{T-2} \beta_{s,i}(u)}{1-\beta_{s,T-1}(u) } .$$ 
\end{lem}

{\it Proof of Lemma \ref{mixingale}.} We use here a limit theorem for $\L^1$-mixingales established in \cite{and88}. Indeed, for $u \in (0,1)$, $s\in \{1,\ldots,T\}$, let 
\begin{equation}\label{eqZ}
Z_{n,t}=\frac 1 { b_n} \,H\Big(\frac{\frac {s+(t-1)T}{nT}-u} {b_n}\Big) \Big ( \big (X^{(n)}_{s+(t-1)T-1}\big )^2- \E \big (X^{(n)}_{s+(t-1)T-1}\big )^2 \Big ).
\end{equation}
Then, set $$c_0(t)=1, \ \ \mbox{ and }\quad c_k(t)=\prod_{i=1}^k a_{t+1-i}\Big(\frac {t+1-i}{nT}\Big),\quad \mbox{ for }k\geq 1,$$ we have:
\begin{eqnarray}\label{lin}
X^{(n)}_{t}=\sum_{k=0}^\infty c_k(t) \, \xi_{t-k}.
\end{eqnarray}
Therefore, with $({\cal F}_{n,t}^{(s)})$ defined in \eqref{F},
\begin{multline*}
\E \big [ Z_{n,t} | {\cal F}_{n,t-m}^{(s)} \big ]=\frac 1 {b_n} \,H\Big(\frac{\frac {s+(j-1)T}{nT}-u} {b_n}\Big) \\
\times\Big \{ \E \Big [ \Big ( \sum_{k=0}^\infty c_k(s+(t-1)T-1) \, \xi_{s+(t-1)T-1 -k} \Big )^2  | ~{\cal F}_{n,t-m}^{(s)}  \Big ] \\
 -\sigma^2 \, \sum_{k=0}^\infty c^2_k(s+(t-1)T-1)   \Big \} \\
  = \frac 1 {b_n} \,H\Big(\frac{\frac {s+(j-1)T}{nT}-u} {b_n}\Big) \Big \{  \Big ( \sum_{k=mT-1}^\infty c_k(s+(t-1)T-1) \, \xi_{s+(t-1)T-1 -k} \Big )^2  \\ 
-\sigma^2 \, \sum_{k=mT-1}^\infty c^2_k(s+(t-1)T-1)   \Big \}.
\end{multline*}
But for any $t \in \N$, we have $|c_k(t)| \leq \alpha^k$ from Assumption (A$(\rho)$). Then, 
\begin{multline*}
\left\| \E \big [ Z_{n,t} ~| ~{\cal F}_{n,t-m}^{(s)} \big ]\right\| _1 \leq  \frac 1 {b_n} \left|H\Big(\frac{\frac {s+(j-1)T}{nT}-u} {b_n}\Big)\right|
\\
\times  \Big \{ \E  \Big [ \Big ( \sum_{k=mT-1}^\infty c_k(s+(t-1)T-1) \, \xi_{s+(t-1)T-1 -k} \Big )^2 \Big ]
 +\sigma^2  \sum_{k=mT-1}^\infty \alpha^{2k} \Big \} \\
 \leq  \frac {2 \,\sigma^2}  {b_n}\left|H\Big(\frac{\frac {s+(j-1)T}{nT}-u} {b_n}\Big)\right| \times \frac {\alpha ^{2mT-2}}{1 -\alpha^2}.
\end{multline*}

Thus, using the notations of Definition 2 in \cite{and88}, it is easy to derive that $(Z_{n,t})$ is a triangular array such that  
$\phi_m=\alpha ^{2mT-2}\|H\|_1 \to 0$ (as $m\to \infty$) since $0\leq \alpha<1$ and: 
$$
\frac 1 n \, \sum_{t=1}^n |c_{nt}| \limiten \frac {2 \,\sigma^2}  {(1 -\alpha^2) }\|H\|_1<\infty, \quad \mbox{with}\quad c_{nt}=\frac {2 \,\sigma^2}  {(1 -\alpha^2) b_n} \,H\Big(\frac{\frac {s+(j-1)T}{nT}-u} {b_n}\Big).
$$ 
As a consequence,
$$\frac 1 n \, \sum_{t=1}^n Z_{n,t} \limiteproban 0,  $$
implies 
$$
  \frac 1 {n  b_n}  \sum_{j=1}^n H\Big(\frac{\frac {s+(j-1)T}{nT}-u} {b_n}\Big) \Big ( \big (X^{(n)}_{s+(j-1)T-1}\big )^2 -\E \big ( \big (X^{(n)}_{s+(j-1)T-1}\big )^2 \big ) \Big ) \limiteproban 0.
$$
Now, we collect the above relations. Lemma \ref{Riemann} and  Proposition \ref{lem1} with  the  $\rho-$regularity of the function $c(v)$,  together conclude the proof.\findem
\begin{lem}\label{Orlicz}
Under the conditions of Theorem \ref{theo1}, with $(Y_{n,i})_ {1\leq i \leq n,\, n\in \N}$ defined in \eqref{Y}, for any $\varepsilon>0$,
\begin{eqnarray}\label{convLindeb}
 \sum_{j=1}^n \E \big ( Y_{n,j}^2 \,\1_{\{ |Y_{n,j}| \geq \varepsilon\}} | \, {\cal F}^{(s)}_{j-1}\big ) \limiteproban 0.
\end{eqnarray}
\end{lem}
{\it Proof of Lemma \ref{Orlicz}.}  Since $\E \xi_0^2  =\sigma^2<\infty$ this is easy to exhibit an increasing sequence $(c_k)_k$ with 
$$
c_0=1,~c_1=2 ~\mbox{and}~c_{k+1} \geq c_k^2,\mbox{ where }
\E \big (\xi_0^2 \,\1_{\{|\xi_0| \geq c_k\}}\big ) \leq \frac 1 {k^3},\quad \mbox{ for all } k \in \N^*.$$
Define  
$g( \cdot)$ as the piecewise affine function such that 
$g(c_k)=k$ for $k \in \N$ and $g(0)=0$. Then the function $\psi$ defined by
$\psi(x)=x^2g(x)$ for $x\geq 0$ satisfies $\psi(0)=0$ and it is a continuous and non-decreasing function (for almost all $x>0$, $\psi'(x)=x^2g'(x)+2xg(x)>0$) and convex function (indeed, for almost all $x>0$, $\psi''(x)=4xg'(x)+2g(x)>0$). Hence, we have:
\begin{multline*}
\sum_{k=1} ^\infty \E \Big ( \xi_0^2 g(|\xi_0|)   \1_{\{c_k\leq |\xi_0| <c_{k+1}\}} \Big )
 \leq   \sum_{k=0} ^\infty \E \Big ( \xi_0^2 g(|\xi_0|)   \1_{\{k\leq g(|\xi_0|) <k+1\}} \Big)\\
 \leq  \sum_{k=1} ^\infty (k+1)  \E \big ( \xi_0^2 \1_{\{c_k\leq |\xi_0|\}} \big) \leq   \sum_{k=1} ^\infty \frac {k+1} {k^3}  <\infty.
\end{multline*}
Therefore, 
\begin{eqnarray}
\label{psi} \E \psi(|\xi_0|) \leq \E \big ( \xi_0^2 g(|\xi_0|)\1_{\{0\leq |\xi_0| <2 \}} \big )+ \sum_{k=0} ^\infty \E \big ( \xi_0^2 g(|\xi_0|)  \1_{\{c_k\leq |\xi_0| \big)<c_{k+1}\}}  
< \infty. 
\end{eqnarray}
The construction of $(c_k)_k$ and the relation  $c_{k+1} \geq c^2_k$  together imply:
\begin{equation}\label{orl2} 
\psi(xy)\le  \psi(x) \psi(y).
\end{equation}
Indeed, this relationship is equivalent to
\begin{eqnarray}\label{gxy}
g(xy)\leq g(x)g(y),\quad \mbox{for any} ~0\leq x\leq y.
\end{eqnarray}
But if $0\leq x \leq 1$ and $y\geq x$, then $xy\leq y$: therefore $g(xy)\leq g(y) \leq g(x)g(y)$ since $g$ is an increasing function and $g(x)\geq 1$ for any $x\geq 0$. Moreover, if $1<x\leq y$, there exists $0\leq k$ and $\lambda \in [0,1]$ such as $y=\lambda c_k+(1-\lambda)c_{k+1}$. But $h:[0,\infty)\to \R^+$ defined by $x  \mapsto h(x)=g(x^2)$ is a convex function since $h''\geq 0$ a.e. As a consequence, 
\begin{multline*}
g(y^2)=h(\lambda c_k+(1-\lambda)c_{k+1})\leq \lambda g(c_k^2)+(1-\lambda) g(c_{k+1}^2)\\
\leq \lambda g(c_{k+1})+(1-\lambda) g(c_{k+2}) \leq \lambda (k+1)+(1-\lambda)(k+2)=k+2-\lambda,
\end{multline*} from the construction of $(c_k)$. Since $g(y)=\lambda g(c_k)+(1-\lambda) g(c_{k+1})=k+1-\lambda$ because $g$ is a piecewise function, we finally obtain $g(y^2)\leq g(y)+1$. We conclude with $g(xy)\leq g(y^2)$ for any $1 \leq x\leq y$ and $g(x)\geq 2$ (since $c_1=2$).  \\
Hence the function $\psi $ is a Orlicz function and $\|\xi_0\|_\psi<\infty$
with
\begin{equation}\label{orl}
\|V\|_\psi=\inf\Big \{z>0;\;\E \psi\Big(\frac{|V|}z\Big)\le1 \Big \},\qquad \mbox{for any random variable $V.$}
\end{equation}
Now Theorem 1.1 in \cite{mal} implies:
\begin{equation}\label{orl3} 
 \|V\|_\psi\le \inf_{z >0}  \frac 1 z   \big(1+\E\big [\psi(z   |V|) \big ]\big)\le 2  \, \|V\|_\psi.
\end{equation}
Therefore $ \|V\|_\psi \le 1+\E\psi(|V|) $, and   $\frac 1 z   \E \psi(z  |V|)\le 2    \|V\|_\psi$ for any $z> 0$ since from convexity  $$\E \psi(|V|) \le \frac{z-1}z    \cdot \E\psi(0)+\frac1 z\cdot   \E \psi(z |V|) \le 2   \|V\|_\psi$$ and $\psi(0)=0$.
\bigskip
Then, from the definition of $(X_t^{(n)})$ and the triangular inequality 
$$
\|X^{(n)}_t\|_\psi\le \alpha    \|X^{(n)}_{t-1}\|_\psi+ \|\xi_t\|_\psi\le  \sum_{j=0}^{t-1} \alpha^j  \|\xi_{t-j}\|_\psi \quad \mbox{for any $t \in \N^*$},
$$
with $0\leq \alpha<1$. Since $\|\xi_{s}\|_\psi=\|\xi_{0}\|_\psi$ for any $s \in \N$, we finally obtain 
$$
\sup_{t \in \N}  \big \{ \|X_t^{(n)}\|_\psi \big \}\le  \frac 1 {1-\alpha}   \|\xi_0\|_\psi <\infty.
$$
Thus \eqref{orl2} implies with the independence of $\xi_t$ and $X^{(n)}_{t-1}$ that:
$$
\E  \psi(|\xi_t   X^{(n)}_{t-1}|)  \le \E \psi(|\xi_t|) \cdot \E  \psi(|X^{(n)}_{t-1}|) .
$$
Now relation \eqref{orl3} with $z=1$ entails
$$
\sup_{t \in \N^*} \big \{ \|\xi_t   X^{(n)} _{t-1}\|_\psi \big \}<\infty.
$$
Thus with $t=s+(j-1)T$ we have from   \eqref{orl3},
$$
\|Y_{n,j}\|_\psi \leq \frac 1 {\sqrt{n b_n} }    \Big |K\Big(\frac{\frac {t}{nT}-u} {b_n}\Big) \Big | \|\xi_t\|_\psi\|X^{(n)}_{t-1}\|_\psi<\infty.
$$
Again using  \eqref{orl2} and with $K_t= K\Big(\frac{\frac {t}{nT}-u} {b_n}\Big)$, 
\begin{eqnarray*}
\E \big ( Y_j^2   \1_{\{|Y_j| \geq \varepsilon\}}\big ) & = &   \frac 1 {nb_n}    \E \big ( (K_t   \xi_t   X^{(n)}_{t-1})^2   \1_{\left\{g(|K_t   \xi_t   X^{(n)}_{t-1}|) \geq  g(\varepsilon   \sqrt{n   b_n})\right\}}\big )  \\
& \leq &  \frac 1 {nb_n}    \E \Big ( (K_t   \xi_t   X^{(n)}_{t-1})^2  \cdot \frac { g(|K_t   \xi_t   X^{(n)}_{t-1}|)}{g(\varepsilon   \sqrt{n   b_n})}    \1_{\left\{g(|K_t   \xi_t   X^{(n)}_{t-1}|) \geq  g(\varepsilon   \sqrt{n   b_n})\right\}}\Big ) \\
& \leq & \frac 1 {\psi(\varepsilon   \sqrt{n   b_n})}    \E \big ( \psi(K_t   \xi_t   X^{(n)}_{t-1})\big ) \\
& \leq & \frac {2   \psi(|K_t |)} {\psi(\varepsilon   \sqrt{n   b_n})}    \sup_{t\in \N^*} \| \xi_t   X^{(n)}_{t-1}\|_\psi. 
\end{eqnarray*}
As a consequence, for any $\varepsilon>0$, 
\begin{multline*}
\E \Big (\sum_{j=1}^n \E \big ( Y_{n,j}^2  \1_{\{ |Y_{n,j}|>\varepsilon\}} |   {\cal F}^{(s)}_{j-1}\big ) \Big ) \\
\leq \frac  {\sup_{t\in \N^*} \| \xi_t   X^{(n)}_{t-1}\|_\psi} {\varepsilon^2 \, g(\varepsilon \, \sqrt{n   b_n})}   \times \frac 1 {n  b_n}  \sum_{j=1}^n  \psi \Big ( \Big |K\Big(\frac{\frac {s+(j-1)T}{nT}-u} {b_n}\Big) \Big | \Big )  \\
 \leq 2\times \frac  {\sup_{t\in \N^*} \| \xi_t   X^{(n)}_{t-1}\|_\psi} {\varepsilon^2 \, g(\varepsilon \, \sqrt{n   b_n})}    \int_\R \psi(|K(x)|)\, dx 
\end{multline*}
if $n$ is large enough, from Lemma \ref{Riemann}. As a consequence, for any $\varepsilon>0$, since $g(\varepsilon\,\sqrt{n   b_n}) \limiten \infty$, then $\E \Big (\sum_{j=1}^n \E \big ( Y_{n,j}^2  \1_{ \{|Y_{n,j}|>\varepsilon\}} |   {\cal F}^{(s)}_{j-1}\big ) \Big ) \limiten 0$. Since $Y_{n,j}^2  \1_{\{ |Y_{n,j}|>\varepsilon\}}$ is a non-negative triangular array, the proof of Lemma \ref{Orlicz} is complete.
\findem

{\it Proof of Theorem \ref{theo1}.}
Using \eqref{AR1}, write
\begin{eqnarray*}\nonumber
\widehat N^{(n)}_s(u)=\frac 1 {nb_n }   \sum_{j\in I_{n,s}}K\Big(\frac{\frac j{nT}-u} {b_n}\Big )   X^{(n)}_{j-1}\Big (a_{s}\Big(\frac j{nT}\Big)X^{(n)}_{j-1}+\xi_{j}\Big)
\end{eqnarray*}
 we decompose it as:
$\widehat N^{(n)}_s(u)=\widetilde N^{(n)}_s(u)+M^{(n)}_{s}(u)$, with
\begin{eqnarray*}
M^{(n)}_{s}(u)&=&\frac 1 {nb_n }  \sum_{j\in I_{n,s}}K\Big(\frac{\frac j{nT}-u} {b_n}\Big) \xi_{j}X^{(n)}_{j-1},\\
\widetilde N^{(n)}_s(u)&=&\frac 1 {nb_n }   \sum_{j\in I_{n,s}}K\Big (\frac{\frac j{nT}-u} {b_n}\Big)a_{s}\big(\frac {j}{nT}\big) (X^{(n)}_{j-1})^{2}
\end{eqnarray*}
Therefore
 we obtain: 
\begin{eqnarray}\label{Decoup} \sqrt{n  b_n}   \big ( \widehat a_s(u)-a_s(u)  \big )&=& \sqrt{n  b_n}    \, \frac {M^{(n)}_{s}(u)}{\widehat D^{(n)}_s(u)} + \frac {J_n} {\widehat D^{(n)}_s(u)},
\end{eqnarray}
with
\begin{eqnarray}
\label{D_n} 
\widehat D^{(n)}_s(u)&=&   \frac1{nb_n }\sum_{j\in I_{n,s}}
K\Big(\frac{\frac j{nT}-u} {b_n}\Big) X_{j-1}^2, \\
\label{J_n}
 J_n\quad\  &=&\frac 1 {\sqrt{ n b_n}}   \sum_{j\in I_{n,s}}K\Big (\frac{\frac j{nT}-u} {b_n}\Big)  (X^{(n)}_{j-1})^{2}   \Big ( a_{s}\big(\frac {j}{nT}\big)-a_s(u) \Big ).
\end{eqnarray}
We are going to derive the consistency of the estimator $\widehat a_s(u)$ of $a_s(u)$, in two parts.
\begin{description}
\item [1/ ] 
We  first prove that $\sqrt{n  b_n}   {M^{(n)}_{s}(u)}\Big/{\widehat D^{(n)}_s(u)} \limiteloin {\cal N} \big ( 0   ,   C \big )$ for some convenient constant  $C>0$. \\
Let $s\in \{1,\ldots,T\}$ and $u\in (0,1)$. For $n \in \N^*$ and $j\in\{1,\ldots,n\}$, we denote
\begin{eqnarray}\label{Y}
Y_{n,j}= \frac 1 {\sqrt{n b_n}}   K\Big(\frac{\frac {s+(j-1)T}{nT}-u} {b_n}\Big) \xi_{s+(j-1)T}X^{(n)}_{s+(j-1)T-1}.
\end{eqnarray}
This is clear that $(Y_{n,j})_{\leq j \leq n,~n\in \N^*}$ is a triangular array of martingale increments with respect to the $\sigma$-algebra ${\cal F}^{(s)}_{t}=\sigma\big ((\xi_{i})_{i\leq s+(t-1)T} \big )$. Indeed $(X^{(n)}_t)_{t\ge0}^{\ }$ is a process, causal with respect to $(\xi_t)_{t\ge0}$. This implies that $\xi_t$ is independent of $\displaystyle (X_i^{(n)})_{i\leq t-1}$ and that $\E(\xi_0)=0$. We are going to use a central limit theorem for triangular arrays of martingale increments, see for example \cite{heyde} and more recently \cite{majorclt}.  
\\
Denote
\begin{eqnarray*}
\sigma_{n,j}^2=  \E \big ( Y_{n,j}^2   | \, {\cal F}^{(s)}_{j-1}\big )  
 = \frac 1 {n b_n}   K^2\Big(\frac{\frac {s+(j-1)T}{nT}-u} {b_n}\Big) \big (X^{(n)}_{s+(j-1)T-1}\big )^2,
\end{eqnarray*}
since $\E(\xi^2_0)=0$. Using Lemma \ref{mixingale}, we obtain:
\begin{eqnarray}\label{convL2}
 \sum_{j=1}^n \sigma_{n,j}^2 \limiteproban \sigma^2  \cdot\frac{ 1+ \sum _{i=0}^{T-2} \beta_{s,i}(u)}{1-\beta_{s,T-1}(u) }     \cdot \int_ \R K^2(x)  dx,
\end{eqnarray}
$\widehat D^{(n)}_s(u)$ is defined from \eqref{D_n} and satisfies
\begin{eqnarray}\label{convL2bis}
\widehat D^{(n)}_s(u) \limiteproban \sigma^2  \, \frac{ 1+ \sum _{i=0}^{T-2} \beta_{s,i}(u)}{1-\beta_{s,T-1}(u) }\equiv\gamma^{(2)}_s(u).
\end{eqnarray}
Moreover, from Lemma \ref{Orlicz}, then for any $\varepsilon>0$,
$$\displaystyle
 \sum_{j=1}^n \E \big ( Y_{n,j}^2  \1_{\{ |Y_{n,j}| \geq \varepsilon\}} | \, {\cal F}^{(s)}_{j-1}\big ) \limiteproban 0
.$$
As a consequence, the conditions of the central limit theorem for triangular arrays of martingale increments, in  \cite{majorclt}), are satisfied and this implies that $\displaystyle \frac {\sum _{j=1}^n Y_{n,j}} {\sqrt {\sum _{j=1}^n \sigma_{n,j}^2} } \limiteloin  {\cal N}\big ( 0  , 1 \big ) $. 
\\
Therefore from Slutsky lemma entails:
\begin{multline}
  \sqrt {nb_n}   \frac {M^{(n)}_{s}(u)}{\widehat D^{(n)}_s(u)} \\
 = \frac {\sum _{j=1}^n Y_{n,j}} {\sqrt {\sum _{j=1}^n \sigma_{n,j}^2} }
\times  \frac {\sqrt { \sum _{j=1}^n \sigma_{n,j}^2} }{\frac 1 {n   b_n}   \sum_{j=1}^n K\Big(\frac{\frac {s+(j-1)T}{nT}-u} {b_n}\Big) \big (X^{(n)}_{s+(j-1)T-1}\big )^2} \\
 \limiteloin  {\cal N}\left ( 0   ,\sigma^2 \,   \frac { 1-\beta_{s,T-1}(u)}{1+ \sum _{i=0}^{T-2} \beta_{s,i}(u) }\,  \int_ \R K^2(x) \, dx \right )\label{convM}
.\end{multline}

\item [2/ ]    The second term $ J_n/\widehat D^{(n)}_s(u)$ in the expansion of $\sqrt{nb_n} \big (\widehat a_s(u)-a_s(u) \big )$ depends on the non-martingale term $J_n$, see \eqref{J_n}, and the consistent term $\widehat D^{(n)}_s(u)$, see \eqref{D_n} and \eqref{convL2bis}. The asymptotic behavior of this second term can be first obtained  following two steps. 
\begin{description}
\item[a. ] A first step consists in establishing an expansion of  $\E J_n$. Using Proposition \ref{lem1} and with $\gamma_s^{(2)} \in {\cal C}^\rho([0,1])$ defined in \eqref{vt2}, we have
\begin{multline*}
\E J_n= \sqrt{ n b_n} \, \frac 1 {n b_n}  \, \sum_{j\in I_{n,s}}K\Big (\frac{\frac j{nT}-u} {b_n}\Big) \\ \times \Big ( \gamma_s^{(2)}(\frac j{nT})+ {\cal O}\big ( \frac 1 n \big )   \Big ) \Big ( a_{s}\big(\frac {j}{nT}\big)-a_s(u) \Big ).
\end{multline*}
Using twice Lemma \ref{Riemann}, with firstly  $c(x)=\gamma_s^{(2)}(x)( a_{s}(x)-a_s(u))$, and secondly $c(x)=( a_{s}(x)-a_s(u))$, we deduce:
\begin{eqnarray}\label{bias}
\big |\E J_n \big  |&\leq & C \,\sqrt{ n b_n} \, \Big ( \frac { A_n } {nb_n}  +b_n^{\rho} \Big )\Big ( 1+ {\cal O}\big ( \frac 1 n \big )   \Big ) .
\end{eqnarray}
As a consequence, if $b_n=o\big ( n^{-1/(1+2\rho)}\big )$, then $\E J_n \limiten 0$.\\ 
In the case  $\rho \in \{1,2\}$, we also obtain from \eqref{equbiais} and with $d_s(v)=(a_s(v)-a_s(u))\gamma_s^{(2)}(v) \in {\cal C}^\rho([0,1])$,
\begin{eqnarray}
\nonumber \E J_n&=&\sqrt{nb_n} \, \Big ({\cal O}\big (\frac {A_n}{nb_n}\big )+
b_n^\rho\, \frac{d_s^{(\rho)}(u)}{\rho!} \int_\R z^\rho K(z)\,dz \,\big (1+o(1)\big )\Big)\\
\nonumber  &=&
\frac{d_s^{(\rho)}(u)}{\rho!}\int_\R z^\rho K(z)\,dz \, \sqrt{nb_n^{2\rho+1}}\, +o(\sqrt{nb_n^{2\rho+1}})+{\cal O}\big (\frac {A_n}{\sqrt{nb_n}}\big )\\
\label{bias2} &=& \left \{ \begin{array}{ll} o(\sqrt{nb_n^{3}})+{\cal O}\big (\frac {A_n}{\sqrt{nb_n}}\big ),&\mbox{if $\rho=1$} ,\\
B_s(u)\,\sqrt{nb_n^{5}} +o(\sqrt{nb_n^{5}})+{\cal O}\big (\frac {A_n}{\sqrt{nb_n}}\big ),&\mbox{if $\rho=2$}.
\end{array} \right .
\end{eqnarray}
with $\displaystyle B_s(u)=\frac{d_s''(u)}{2}\,\int_\R z^2 K(z)\,dz $.
\item[b. ] Now we are going to prove a first consistency result for $ J_n/\widehat D^{(n)}_s(u)$ using the Markov Inequality. Indeed,
\begin{eqnarray*}
\E  |  J_n| &  \leq & \frac 1 {  \sqrt { n b_n}} \sum_{j\in I_{n,s}} \Big | K\Big (\frac{\frac j{nT}-u} {b_n}\Big) \Big | \big | a_{s}\big(\frac {j}{nT}\big)-a_s(u)  \big |   \v (X^{(n)}_{j-1})  \\
& \leq & \frac {C} {  \sqrt { n b_n}}\sum_{j\in I_{n,s}} \Big | K\Big (\frac{\frac j{nT}-u} {b_n}\Big) \Big | \big | a_{s}\big(\frac {j}{nT}\big)-a_s(u)  \big |. 
\end{eqnarray*}
Now using Lemma \ref{Riemann} with $c(v)=| a_{s}(v)-a_s(u) |$ which also belongs in ${\cal C}^\rho([0,1])$ (this is clear if $\rho<1$ and, for $\rho=1$ the Lipschitz property of $z\mapsto |z|$ allows to conclude), and  $c(u)=0$, we derive:
\begin{eqnarray}\label{J1}
\E  |  J_n|  \leq  \sqrt { n b_n}  \Big (\frac {A_n }{nb_n} +b_n^{\rho\wedge 1} \Big ). 
\end{eqnarray} 
\end{description}
Therefore, if $b_n=o\big (n^{-\frac 1{1+2(\rho\wedge 1}} \big )$, then $\E J_n \limiten 0$ and $\E  |  J_n|  \limiten 0$, implying from Markov Inequality, $J_n \limiteproban 0$. Finally, since \eqref{convL2bis} establishes the consistency of $\widehat D^{(n)}_s(u)$, from Slutsky lemma, we deduce 
\begin{equation}\label{second}
\frac {J_n} {\widehat D^{(n)}_s(u)} \limiteproban 0.
\end{equation} 
As a consequence, the proof of the Theorem results by  using the decomposition \eqref{Decoup}, the consistency results \eqref{convM} and \eqref{second}. \findem 
\end{description}
{\it Proof of Theorem \ref{theo2}.} We restrict to the case $\rho\in (1,2]$. 
\begin{description}
\item [a. ] \underline{Case  $\E \big ( \xi_0^4 \big )<\infty$}.  \\
Denote again $\displaystyle K_t= K\Big(\frac{\frac {t}{nT}-u} {b_n}\Big)$, for $t \in \Z$. 
First remark that the symmetry assumption on $\xi_0$'s distribution implies $\E \big (\xi_0\big )=\E \big (\xi_0^3\big )=0$. 
\begin{eqnarray*}
 \v  (  J_n  ) &= &\frac  1 {nb_n}  \, \sum_{t\in I_{n,s}}\sum_{t'\in I_{n,s}} K_t K_{t'} \, \cov(X_t^2, X _{t'}^2)\times \\
&& \qquad \qquad \qquad\times \Big( a_{s}\big(\frac {t}{nT}\big)-a_s(u)  \Big ) \Big ( a_{s}\big(\frac {t}{nT}\big)-a_s(u)  \Big )\\
 & = &\frac  1 {nb_n} \,   \sum_{(t,t')\in L_{n,s,\alpha} } K_t K_{t'} \, \cov(X_t^2, X _{t'}^2)\times \\
&& \qquad \qquad \qquad\times
 \Big  ( a_{s}\big(\frac {t}{nT}\big)-a_s(u)  \Big )   \Big( a_{s}\big(\frac {t'}{nT}\big)-a_s(u)  \Big )   \\
& +& \frac  1 {nb_n} \,   \sum_{(t,t')\in I_{n,s}^2 \setminus L_{n,s,\alpha}}\!\!\!\!K_t K_{t'}  \cov(X_t^2, X _{t'}^2)\times \\
&& \qquad \qquad \qquad\times
\Big  ( a_{s}\big(\frac {t}{nT}\big)-a_s(u)  \Big )   \Big ( a_{s}\big(\frac {t'}{nT}\big)-a_s(u) \Big )  
\end{eqnarray*}
with  $L_{n,s,\alpha}=\big \{(t,t')\in I^2_{n,s},~ \, |t-t'| \leq \frac{\log n}{\log\alpha} \big \}$. \\
Firstly, consider the first left side term of the last inequality.  If $t \in I_{n,s}$  then Proposition \ref{lem1} entails $\v(X_t^2)=\gamma_s^{(2)}(t/(nT)) + {\cal O}(1/n)$ for an adequate function $\gamma_s^{(2)} \in {\cal C}^\rho([0,1])$.\\ Hence we also have $\v(X_{t}^2)=\gamma_s^{(2)}(t/(nT)) + {\cal O}(\log(n)/n)$.\\ Here the fact that $(z\mapsto z^2)$ is a function in $ {\cal C}^\rho$, implies that the function defined from $b(v)=\big ( a_{s}(v)-a_s(u)  \big ) ^2$ is in ${\cal C}^\rho([0,1])$ too, and again $b(u)=0$ and $\int xH^2(x)dx=0$.\\
 Therefore, we  use Lemma \ref{Riemann} to derive:

\begin{multline*}  \nonumber
 \sum_{t,t'\in L_{n,s,\alpha}} K_t K_{t'}  \cov(X_t^2, X _{t'}^2) \Big  ( a_{s}\big(\frac {t}{nT}\big)-a_s(u)  \Big ) \Big ( a_{s}\big(\frac {t'}{nT}\big)-a_s(u)  \Big ) 
 \\ \nonumber
=   \sum_{t,t'\in L_{n,s,\alpha}} K_t K_{t'}  \prod_{i=1}^{|t-t'|} a^2_{s+i} \big (\frac {t} {nT}  \big )   \Big ( \gamma_{s}^{(4)}
\big (\frac {t+i}{nT} \big ) +{\cal O}\big (\frac 1 n \big ) \Big )
\\ \nonumber
\qquad\qquad\qquad  \times \Big  ( a_{s}\big(\frac {t}{nT}\big)-a_s(u)  \Big ) \Big  ( a_{s}\big(\frac {t'}{nT}\big)-a_s(u)  \Big )
\\ \nonumber
=  \sum_{j=0}^{\frac{\log n}{\log\alpha}}  \sum_{t \in I_{n,s}} K_t\Big (K_t+{\cal O}\big (\frac {\log n} {nb_n^2} \big ) \Big )\prod_{i=1}^{j}  a^2_{s+i} \big (\frac {t} {nT}  \big )  \Big ( \gamma_{s}^{(4)}\big (\frac {t}{nT} \big ) +{\cal O}\big (\frac {\log n} n \big ) \Big )
 \\
\qquad\qquad\qquad \times \Big  ( a_{s}\big(\frac {t}{nT}\big)-a_s(u)  \Big )   \Big( a_{s}\big(\frac {t'}{nT}\big)-a_s(u)  \Big )\\ \nonumber
 \leq  \sum_{j=0}^{\frac{\log n}{\log\alpha}} \alpha^{2j}   \Big (  \sum_{t \in I_{n,s}} K^2_t \prod_{i=1}^{j}  \gamma_{s}^{(4)}\big (\frac {t}{nT} \big ) \times \Big  ( a_{s}\big(\frac {t}{nT}\big)-a_s(u)  \Big ) ^2+ {\cal O}\big (\frac {\log n} {n b_n^2} \big )   \Big ) \Big ) \nonumber \\
 \leq  2  \sum_{j=0}^{\infty}  \alpha^{2j}    \sum_{t \in I_{n,s}} K^2_t \prod_{i=1}^{j}   \gamma_{s}^{(4)}\big (\frac {t}{nT} \big ) \times \Big  ( a_{s}\big(\frac {t}{nT}\big)-a_s(u)  \Big ) ^2\nonumber \\
\qquad \leq  2 nb_n   \sum_{j=0}^{\infty}  \alpha^{2j}    \frac 1 {nb_n} \sum_{t \in I_{n,s}} K^2 \Big(\frac{\frac {t}{nT}-u} {b_n}\Big)  g_j \big (\frac {t}{nT}\big) ,\end{multline*}
with $\displaystyle g_j(x) = \big  ( a_{s}(x)-a_s(u)  \big ) ^2\prod_{i=1}^{j}  \big ( \gamma_{s}^{(4)} (x )$,
since for $n$ large enough the above expression satisfies  $\big |{\cal O}\big (\frac {\log n} {n b_n^2} \big ) \big |\leq 1$. Using  Lemma \ref{Riemann}, with functions $H=K^2$ and $c=g_j$ with $g_j \in {\cal C}^\rho ([0,1])$ (quote that $\max _{i\le j} \big ( \| g_i \| \vee  \Lip (g_i) \big )={\cal O }(j)$), we finally obtain:
\begin{multline}\label{T1} 
\Big |\frac 1 {nb_n}  \sum_{t,t'\in L_{n,s,\alpha}} K_t K_{t'}  \cov(X_t^2, X _{t'}^2) \big  ( a_{s}\big(\frac {t}{nT}\big)-a_s(u)  \big )  ( a_{s}\big(\frac {t'}{nT}\big)-a_s(u)  \big ) \Big | \\
\leq C    \Big (\frac {A_n }{nb_n}  +b_n^{\rho} \Big ) .
\end{multline}
Secondly, from Proposition \ref{lem1}, for $t,t' \in I_{n,s}^2\setminus L_{n,s,\alpha}$, we have $$ |\cov (X_t^2,X_{t'}^2) |\leq C \, \alpha ^{2\,|t-t'|} \leq  \frac C{n^2}.$$ Thus,
\begin{multline} \label{T2}\!\!\!\!\!\!\!\!\!
\Big |\frac 1 { n b_n}    \sum_{t,t'\in I^2_{n,s} \setminus L_{n,s,\alpha}}\! \!\!\! K_t K_{t'}   \cov(X_t^2, X _{t'}^2)  \Big  ( a_{s}\big(\frac {t}{nT}\big)-a_s(u)  \Big ) \Big ( a_{s}\big(\frac {t'}{nT}\big)-a_s(u)  \Big )   \Big | \\
\leq \frac {nb_n} {n^2} \, \Big ( \frac 1 {n b_n} \,  \sum_{t\in I_{n,s}}K_t \big  ( a_{s}\big(\frac {t}{nT}\big)-a_s(u)  \big )  \Big ) ^2 \leq C \, \frac {b_n } n \, \Big (\frac { A_n}{nb_n}  +b_n^{\rho} \Big ),
\end{multline}
from Lemma \ref{Riemann}. Then, \eqref{T1} and \eqref{T2} provide 
\begin{equation}\label{var}
\v \big (  J_n \big )  \le C  \,  \Big (\frac { A_n}{nb_n}  +b_n^{\rho} \Big )
\end{equation}
implying $\v \big (  J_n \big ) \limiten 0$ for any $(b_n)$ such as $$\max(b_n,A_n(n\, b_n)^{-1} ) \limiten 0.$$ 

\item[b. ] \underline{Case $\E \big ( |\xi_0|^\beta \big )<\infty$, for some $\beta \in [2,4]$.} \\
From its expression given in \eqref{J_n}, $J_n$ is a quadratic form of $(X_t)$ and therefore, as $X_t$ is a linear process with innovations $(\xi_t)$, $J_n$ is also a quadratic form of $(\xi_t)$. As a consequence, the fourth order moment can be injected such as there exists a sequence $z_n\downarrow0$ (as $ n\uparrow\infty$) satisfying:
\begin{equation}\label{J2}
\v (  J_n )\le z_n \, \big (\E(\xi_0^4)\vee 1 \big )=z_n\, \big (\mu_4\vee 1\big ), \mbox{ and } z_n={\cal O}\big (\frac{A_n}{nb_n}+b_n^{\rho}\big).
\end{equation}  
Now, assume only that $\E (\xi_0^2)<\infty$. The innovations $(\xi_t)$ can be truncated at level $M$, and  write 
$$
\xi_{t,M}=\xi_t\, \1 _{|\xi_t|\leq M }\qquad\mbox{for any $t\in \N$}.
$$ 
Note that the symmetry assumption entails  $\E (\xi_{j,M})=0$. Define also
Define also
\begin{eqnarray*}
X_{t,M}^{(n)} &=&  a_{t}\big(\frac t{nT}\big)X^{(n)}_{t-1,M}+\xi_{t,M},\qquad 1\le t\le nT, ~n\in \N \\
 \mbox{and} \quad J_{n,M}& =& \frac 1 {\sqrt{ n b_n}}   \sum_{j\in I_{n,s}}K\Big (\frac{\frac j{nT}-u} {b_n}\Big)  (X^{(n)}_{j-1})^{2}   \Big ( a_{s}\big(\frac {j}{nT}\big)-a_s(u) \Big ).
\end{eqnarray*}  
A consequence of \eqref{J2} is:
\begin{equation}\label{varM}
\v (  J_{n,M} )\le z_n \, \E(\xi_{0,M}^4) \leq z_n \,M^2 h(M),
\end{equation}
with $h(M)= \E \big (|\xi_{0}|^2\1_{\{|\xi_{0}|>M\}} \big )$ which satisfies $\lim_{M\to \infty}h(M)=0$.\\ 
Moreover,
\begin{multline}
 \big |J_n- J_{n,M}\big| 
   = \frac 1 {\sqrt{ n b_n}}    \sum_{j\in I_{n,s}}K\Big (\frac{\frac j{nT}-u} {b_n}\Big) 
   \\  \times\big |(X^{(n)}_{j-1})^{2} -(X^{(n)}_{j-1,M})^{2} \big |   \big | a_{s}\big(\frac {j}{nT}\big)-a_s(u) \big |. 
  \label{troncJ} \end{multline}
But
\begin{eqnarray}\nonumber
X^{(n)}_{j-1,M}-X^{(n)}_{j-1}\;&=&a_{j-1}\big (\frac {j-1}{nT} \big )\big (X^{(n)}_{j-2,M}-X^{(n)}_{j-2}\big )+\big (\xi_{j-1,M}-\xi_{j-1} \big), \\
\label{xm}
 \big |X^{(n)}_{j-1,M}-X^{(n)}_{j-1} \big |&\le&\alpha \, \big  |X^{(n)}_{j-2,M}-X^{(n)}_{j-2} \big |+|\xi_{j-1}|\1_{\{|\xi_{j-1}|>M\}}. 
\end{eqnarray}
We first remark from Proposition \ref{lem1} that $\E (X^{(n)}_{j-1})^{2}+\E (X^{(n)}_{j-1,M})^{2}\le c$ for some constant $c>0$. Hence, Cauchy-Schwartz Inequality shows that, for each $j$:
\begin{equation} \label{deltajM}
\E \big ( \big|(X^{(n)}_{j-1})^{2} -(X^{(n)}_{j-1,M})^{2} \big| \big )\le \sqrt {c\delta_{j-1,M}},\end{equation}
with $\delta_{j-1,M}=\E \big (|X^{(n)}_{j-1} -X^{(n)}_{j-1,M}|^{2}\big )$. 
\\
We are going to bound $\delta_{j-1,M}$. A first simple  bound is clearly $\delta_{j-1,M}\le 2\, c$ and we use it together with \eqref{xm}, and Cauchy-Schwartz inequality in order to derive
\begin{eqnarray*}
\delta_{j-1,M}&\le& \alpha^2 \delta_{j-2,M}+2\alpha  \, \E \big ( |X^{(n)}_{j-2,M}-X^{(n)}_{j-2}| |\xi_{j-1}|\1_{\{|\xi_{j-1}|>M\}}| \big )\\
&&\qquad+ \E \big (|\xi_{j-1}|^2\1_{\{|\xi_{j-1}|>M\}}\big )
\\
&\le& \alpha^2 \delta_{j-2,M}+2\alpha \,\sqrt{2c} \, \sqrt{\E |\xi_{j-1}|^2\1_{\{|\xi_{j-1}|>M\}}}+\E |\xi_{j-1}|^2\1_{\{|\xi_{j-1}|>M\}}
\\
&\le&  \alpha^2 \delta_{j-2,M}+H(M) \quad (\mbox{with}\quad H(M)=2\alpha \,\sqrt{2c} \,  \sqrt{h(M)}+h(M))
\\
&\le&
  \alpha^4 \delta_{j-3,M}+ (1+\alpha^2)H(M)
\\
 &\le&\cdots
 \\
&\le &  \alpha^{2(j-1)} \delta_{0,M}+ (1+\cdots+\alpha^{2(j-2)})H(M)\\
&\le&     \frac2{1-\alpha^2}H(M)
\end{eqnarray*}
since $ \delta_{0,M}\le h(M) \leq H(M)$.\\ Now, from \eqref{deltajM}, we obtain for $M$ large enough:
\begin{equation}\label{eta}
\E  \big|(X^{(n)}_{j-1})^{2} -(X^{(n)}_{j-1,M})^{2} \big|\le  \sqrt{ \frac{2c}{1-\alpha^2}}\,\sqrt{H(M)} \leq C\, h^{1/4}(M)
\end{equation}
with  $C>0$ and always with $h(M)= \E \big (|\xi_{0}|^2\1_{\{|\xi_{0}|>M\}} \big )$. Now a careful use of \eqref{J1} and \eqref{troncJ} entails:
\begin{equation}\label{first}
\E  | J_n- J_{n,M}| \le C\,  \sqrt { n b_n}  \Big (\frac {A_n}{nb_n} +b_n \Big ) \, h^{1/4}(M)
\end{equation}
since $x \to |a(x)-a(u)|$ is a ${\cal C}^1$ function (in the above defined sense). Finally, using Cauchy-Schwartz inequality in \eqref{varM}, we obtain for $M$ large enough,
\begin{eqnarray}
\nonumber \E | J_n| &\le &\E  | J_n- J_{n,M}| + \sqrt{\v ( J_{n,M} )} \\
 \nonumber &\leq & C\, \Big ( \sqrt { n b_n}  \Big (\frac {A_n}{nb_n} +b_n \Big ) \, h^{1/4}(M) + \big (\frac{A_n}{nb_n}+b_n^\rho\big)^{1/2} \,M\, h^{1/2}(M)  \Big ) \\
\label{beta}&\leq & C\, \Big ( \sqrt { n b_n^3} \, h^{1/4}(M) + b_n^{\rho/2} \,M\, h^{1/2}(M)  \Big ) 
\end{eqnarray}
assuming $A_n/nb_n =o(b_n^{\rho/2})$ {\it i.e.} $(n/A_n)^{-2/(2+\rho)}=o(b_n)$ (and note that $-2/(2+\rho) \leq 1/(1+2\rho)$). \\
Now, if $\E \big ( |\xi_0|^{\beta}  \big )<\infty$ with $\beta \in (2,4]$, then using H\"older and Markov Inequalities, there exists $C_\beta>0$ such as 
$$
h(M)= \E \big (|\xi_{0}|^2\1_{\{|\xi_{0}|>M\}} \big ) \leq C_\beta \, M ^{2-\beta}. 
$$ 
Since  here $b_n=o\big (n ^{-1/(1+2\rho)} \big )$, does not  yields  the minimax rates, we deduce that 
$$
\left \{ \begin{array}{lcl} 
\sqrt { n b_n^3} \, h^{1/4}(M) \limiten 0 & \mbox{when} & M^{1+2\rho}\geq n^{(4\rho-4)/(\beta-2)} \\
b_n^{\rho/2} \,M\, h^{1/2}(M)  \limiten 0& \mbox{when} &  M^{1+2\rho}\leq n^{\rho/(4-\beta)}
\end{array} \right . .
$$
Thus, from inequality \eqref{beta}, we deduce that the optimal choice is obtained when  $$\frac{4\rho-4}{\beta-2}=\frac\rho{4-\beta},\ \mbox{ which entails }\beta=4-2\cdot \frac\rho{5\rho-4}.$$ 
\item[d.] {\underline{Case $\rho=2$.}}\\
 The expression of the non-central limit for the case of optimal window widths and the expansion of the bias \eqref{bias2} now the asymptotic expression for \eqref{second}  yields the proposed non-centred Gaussian limit, see Remark \ref{remJM}. The same truncation step as above is also needed.
\end{description}
The proof is now complete.
\findem
\begin{Rem}\label{remJM}
Using the previous bound \eqref{bias} of \ $\E J_n$ and Bienaym\'e-Tchebychev inequality, we deduce that if $b_n=o \big ( n^{-1/(1+2\rho)} \big )$ then $ J_n  \limiteproban 0$.\\
Moreover, if $\rho=2$ and $b_n=c \,  n^{-1/5} $, using the expansion \eqref{bias2} of $\E J_n$ and again Bienaym\'e-Chebychev inequality,  then $J_n \limiteproban B_s(u)\, c^{5/2}$.\\ Therefore with the consistency result \eqref{convL2bis}, for any $u \in (0,1)$ and $s \in \{1,\ldots,T\}$,
$$
\frac {J_n} {\widehat D^{(n)}_s(u)} \limiteproban B_s(u) \, \frac{ c^{\frac52}}{\sigma^2 } \frac{ 1+ \sum _{i=0}^{T-2} \beta_{s,i}(u)}{1-\beta_{s,T-1}(u) }. 
$$
\end {Rem}
\begin{Rem}\label{mom3}
For the general case with maybe $\xi_0$ non symmetric and $\E \xi_0=0$,  the item 3. of Proposition \ref{lem1} needs some improvements. Denote  $w^{(k)}_t=\E(X^k_t)$ for $k=1,3$, then  $w^{(4)}_t=w_t$  and  $w^{(2)}_t=v_t$, then \eqref{w} turns to be written
\begin{equation}\label{w3}
w_t  = q_t w_{t-1} +4\E A_t\E X_{t-1}\mu_3 + 6\sigma^2  v_t + \mu_4 \leq \alpha^4  w_{t-1} + r(t),
\end{equation}
as previously $\sup_{t} w_t <\infty$. \\
We  need to derive suitable equivalents of   $w^{(k)}_t$ if $k=1$. Firstly $$w^{(1)}_{t}=\E A_tw^{(1)}_{t-1}=\cdots=E A_t\cdots\E A_1\cdot \E X_0,$$ and in fact this term is negligible and the proof of Proposition \ref{lem1} and Lemma 3. remains unchanged.
\\
In this case the proof of the  above point {\bf 2/\,c.} needs a simple improvement and
$$\xi_{j,M}=\xi_j\wedge M\vee(-M)-\E\left(\xi_j\wedge M\vee(-M)\right).$$  
In this truncated setting, inequality \eqref{xm}   writes:
\begin{multline*} |X^{(n)}_{j-1,M}-X^{(n)}_{j-1}|\le\alpha |X^{(n)}_{j-2,M}-X^{(n)}_{j-2}|\\
+|\xi_{j-1}|\1_{\{|\xi_{j-1}|>M\}}+ \E(|\xi_{j-1}|\1_{\{|\xi_{j-1}|>M\}})
\end{multline*}
so that the end of the proof is unchanged by only setting  $ C={2c\E\xi_0^2}/{(1-\alpha)}$.

\end{Rem}
\begin{Rem}\label{xm1}
 Secondly, in case we even omit the condition $\E\xi_0=0$ one needs to also express an asymptotic expansion for $w^{(3)}_{t}= \E A_t^3w^{(3)}_{t-1}+3\E A_tw^{(1)}_{t-1}\sigma^2+\mu_3 \sim  \E A_t^3w^{(3)}_{t-1}+\mu_3 $; an analogue expansion to Proposition \ref{lem1} and Lemma 3.  may thus be derived.
Namely
$\displaystyle
w^{(3)}_{t}=\gamma_s^{(3)}(\frac t{nT})+ {\cal O}\big ( \frac 1 n \big ),
$,
with
\begin{eqnarray*}
  \gamma_s^{(3)}(v)&=&\mu_3  \cdot \frac{ 1+ \sum _{i=0}^{T-2} \zeta_{s,i}(v)}{1-\zeta_{s,T-1}(v) } ,\\
 \zeta_{t,i}(v)&=&\prod _{j=0}^{i-1} a^3_{t-j}(v)  \leq \alpha^{3i}<1 ,\qquad \mbox{for} \quad 1 \leq i \leq T, \quad v\in(0,1).
\end{eqnarray*}
Then the expression of the equivalent of $w_t$ is also adequately transformed up to the above relations.
\end{Rem}

\paragraph{Aknowledgement.} 
This work has been developed within the ``MME-DII centre of excellence'' (ANR-11-LABEX-0023-01)  and with the help of PAI- CONICYT MEC Nr. 80170072. \\
The authors thank the referees for their fruitful comments and suggestions,  which notably improved the quality of the paper. The second author wishes to thank Rainer Dahlhaus for many interesting discussions. As well, numerous discussions with Karine Bertin were extremely useful.



\end{document}